\documentclass[11pt]{article}

\usepackage{amssymb}
\usepackage{amsmath,amsfonts,mathtools}
\usepackage{amsthm}
\usepackage{latexsym}
\usepackage{mathrsfs}
\usepackage{color}
\usepackage{float}
\usepackage{enumitem}
\usepackage{algorithm,algorithmic}
\usepackage{graphicx}
\usepackage{subcaption}
\usepackage{booktabs}
 \allowdisplaybreaks[3]
 
\usepackage{bbm}
\usepackage{microtype}
\usepackage[a4paper]{geometry}
\usepackage{pifont}
\usepackage{latexsym}
\usepackage{url}

\usepackage[a4paper]{geometry}
\newcommand{\beq}{\begin{equation}}
\newcommand{\eeq}{\end{equation}}
\newcommand{\beqa}{\begin{eqnarray}}
\newcommand{\eeqa}{\end{eqnarray}}
\newcommand{\beqas}{\begin{eqnarray*}}
\newcommand{\eeqas}{\end{eqnarray*}}
\newcommand{\ba}{\begin{array}}
\newcommand{\ea}{\end{array}}
\newcommand{\bi}{\begin{itemize}}
\newcommand{\ei}{\end{itemize}}
\newcommand{\nn}{\nonumber}

\newcommand{\cK}{{\mathcal K}}
\newcommand{\cL}{{\mathcal L}}
\newcommand{\cQ}{{\mathcal Q}}

\newcommand{\cN}{{\mathcal N}}

\newcommand{\dom}{\mathrm{dom}}
\newcommand{\dist}{\mathrm{dist}}
\newcommand{\argmin}{\arg\min}

\newcommand{\inte}{\mathrm{relint}}

\newcommand{\ceil}[1]{\left\lceil #1 \right\rceil}

\newtheorem{lemma}{Lemma}
\newtheorem{thm}{Theorem}

\newtheorem{defi}{Definition}

\newtheorem{assumption}{Assumption}

\newtheorem{rem}{Remark}

\newcounter{spb}
\def\bepsilon{{\bar \epsilon}}
\def\bK{{\bar K}}

\def\bN{{\bar N}}

\def\cN{{\cal N}}
\def\cO{{\cal O}}

\def\tf{{\tilde f}}
\def\tF{{\tilde F}}

\def\tu{{\tilde u}}

\def\tx{{\tilde x}}
\def\tz{{\tilde z}}

\def\tr0{{\tilde r_0}}

\def\cK{{\mathcal{K}}}

\def\rr{{\mathbb{R}}}

\def\bbE{{\mathbb{E}}}

\def\cB{{\mathscr{B}}}

\def\barf{{\bar f}}
\def\barF{{\bar F}}

\def\inte{\text{int}}
\def\xint{x_{\rm int}}
\def\hx{{\hat x}}

\title{An inexact stochastic variance-reduced gradient method for constrained convex finite-sum optimization }
 \author{
 Zhaosong Lu
 \thanks{
 Department of Industrial and Systems Engineering, University of Minnesota, USA (email: {\tt zhaosong@umn.edu}, {\tt xiao0414@umn.edu}). This work was partially supported by the Office of Naval Research under Award N00014-24-1-2702,  the Air Force Office of Scientific Research under Award FA9550-24-1-0343, and the National Science Foundation under Awards 2211491 and 2435911.
 }
 \and
 Yifeng Xiao
 \footnotemark[1]
 }
 \date{August 5, 2025 (Revised: October 7, 2026)}

\begin{document}
\maketitle

\begin{abstract}
In this paper, we study a class of constrained convex finite-sum optimization problems. In many practical applications, all or some of the constraints are hard constraints that represent imperative requirements. Consequently, any desirable approximate solution must (nearly) satisfy these constraints. Motivated by this, we focus on computing an $\epsilon$-\emph{surely feasible stochastic optimal} ($\epsilon$-SFSO) solution, in which the constraint violation is deterministically bounded by $\epsilon$, and the expected optimality gap is at most $\epsilon$. While a recently proposed stochastic variance-reduced first-order penalty method can compute an $\epsilon$-SFSO solution, there remains a significant gap between its first-order oracle complexity for computing an $\epsilon$-SFSO solution and that for computing an $\epsilon$-\emph{expected feasible stochastic optimal} solution, where both the expected constraint violation and the expected optimality gap are bounded by $\epsilon$. To narrow this gap, we propose an inexact stochastic variance-reduced gradient (SVRG) method that significantly improves the previously best-known first-order oracle complexity for computing an $\epsilon$-SFSO solution. Furthermore, when the objective consists of a single function, the proposed method reduces to a novel deterministic first-order method for solving standard convex conic optimization. It achieves a significantly lower gradient evaluation complexity with respect to the objective function compared to existing methods, and can outperform them when computing constraint gradients { and projections onto the associated cone} is cheaper than evaluating the objective gradient. 
{ Preliminary numerical results are presented to illustrate the performance of the proposed method.}
\end{abstract}

 \noindent{\bf Keywords:} convex conic optimization, finite-sum optimization, accelerated gradient method, variance reduction
 \medskip

 \noindent{\bf Mathematics Subject Classification:}  90C25, 90C35, 90C30, 65K05

\section{Introduction} \label{sec:introduction}

In this paper, we consider constrained convex finite-sum optimization  problems of the form 
\beq\label{prob1}
\begin{aligned}
F^*=\min_{x} &\  \, \Big\{F(x):=\frac{1}{s}\sum_{i=1}^s f_i(x) +\psi(x)\Big\}\\
\mbox{s.t.}&\ \, c(x)\in -\cK,
\end{aligned}
\eeq
where $\cK\subseteq\rr^m$ is a closed convex cone, $f_i$ is a continuously differentiable convex function on $\rr^n$ for all $i=1,2,\dots,s$, $\psi: \rr^n \rightarrow \rr \cup \{+\infty\}$ is a proper closed convex function with an exactly evaluable proximal operator and domain denoted by $X$, and $c = (c_1, \ldots, c_m)$ is a  Lipschitz smooth and $\cK$-convex mapping on $X$, that is,
\beq \label{K-convex}
    \alpha c(x)+(1-\alpha)c(y)-c(\alpha x+(1-\alpha)y)\in\cK \qquad\forall x,y\in X, \alpha\in[0,1].
\eeq

Problem~\eqref{prob1} includes a broad class of constrained optimization problems as special cases. For instance, when $\mathcal{K} = \mathbb{R}_{+}^{m_1} \times \{0\}^{m_2}$ for some $m_1, m_2 \in \mathbb{N}$, and $c = (c_1,\ldots,c_{m_1+m_2})$ with $c_j$ convex for $1 \leq j \leq m_1$ and affine for $m_1 < j \leq m_1 + m_2$, problem~\eqref{prob1} reduces to the standard constrained convex finite-sum optimization problem:
\[
\min_{x} \left\{ \frac{1}{s}\sum_{i=1}^s f_i(x) + \psi(x):\ c_j(x) \leq 0,\ 1 \leq j \leq m_1;\ c_j(x) = 0,\ m_1 <j \leq m_1 + m_2 \right\}.
\]

Moreover, when $c = 0$, problem \eqref{prob1} reduces to the unconstrained finite-sum optimization problem:
\beq \label{fun-f}
\min_{x} \Big\{\frac{1}{s}\sum_{i=1}^s f_i(x) +\psi(x)\Big\}.
\eeq 
This class of problems has been extensively studied in the literature (see, e.g., \cite{allen2018katyusha,defazio2014saga,hendrikx2021optimal,lan2019unified, schmidt2017minimizing,shalev2013stochastic,song2020variance,xiao2014proximal,zhou2019lower}). Notably, several stochastic gradient methods (e.g., \cite{defazio2014saga,schmidt2017minimizing,shalev2013stochastic,xiao2014proximal}) employ variance reduction techniques and achieve $\widetilde{\mathcal{O}}(s \log s+s^{1/2}\epsilon^{-1})$ complexity in terms of the number of $\nabla f_i$ evaluations required to compute an $\epsilon$-\textit{stochastic optimal solution}, that is, a point $x$ satisfying $\mathbb{E}[F(x) - F^*] \leq \epsilon$.\footnote{$\widetilde{\mathcal{O}}(\cdot)$ hides logarithmic factors { of the tolerance  parameter $\epsilon$} in $\mathcal{O}(\cdot)$.} More recently, optimal incremental gradient methods \cite{lan2018optimal} and accelerated stochastic variance-reduced gradient (SVRG) methods \cite{allen2018katyusha,lan2019unified} have been developed, achieving $\widetilde{\mathcal{O}}(s \log s+s^{1/2}\epsilon^{-1/2})$ or $\mathcal{O}(s \log s+s^{1/2}\epsilon^{-1/2})$ complexity for computing such solutions.

In addition, when $c \neq 0$, $\mathcal{K}$ is the nonnegative orthant, and $\psi$ is the indicator function of a simple convex compact set, a level-set method was proposed in \cite{lin2018level} for solving problem~\eqref{prob1}. This method requires $\mathcal{O}\big(sn+(s+n)\epsilon^{-2}\big)$ evaluations of $\nabla f_i$ and $\nabla c$ to compute an $\epsilon$-\emph{expected feasible stochastic optimal} ($\epsilon$-EFSO) solution of~\eqref{prob1}, that is, a stochastic point $x$ satisfying
\begin{equation*} 
\mathbb{E}[\dist(c(x), -\mathcal{K})] \leq \epsilon, \quad \mathbb{E}[F(x) - F^*] \leq \epsilon,
\end{equation*}
where $\dist(\cdot, \cdot)$ denotes the Euclidean distance between a point and a set.

In many applications,  all or some of the constraints in problem \eqref{prob1} are hard constraints that represent imperative requirements. Consequently, any desirable approximate solution must (nearly) satisfy these constraints. Recently, for the case where $\cK$ is the nonnegative orthant, a stochastic variance-reduced first-order penalty (SVRFP) method was proposed in \cite{lu2025first} for computing an $\epsilon$-\emph{surely feasible stochastic optimal} ($\epsilon$-SFSO) solution to problem \eqref{prob1}, that is, a stochastic point $x$ satisfying
\[
\dist\big(c(x),-\cK\big) \leq \epsilon, \quad \bbE[F(x)-F^*] \leq \epsilon,
\]
which is generally a stronger notion than an $\epsilon$-EFSO solution. The SVRFP method applies \emph{only one iteration} of a modified variance-reduced accelerated gradient scheme to a sequence of associated quadratic penalty subproblems:
\beq \label{penalty-prob}
\min_{x}\Big\{\frac{1}{s}\sum_{i=1}^sf_i(x)+\frac{\rho_k}{2}\dist^2\big(c(x),-\cK\big)+\psi(x)\Big\}
\eeq
for a suitably chosen sequence of penalty parameters $\{\rho_k\}$. With an appropriate choice of $\{\rho_k\}$,  this method can compute an $\epsilon$-SFSO solution using $\mathcal{O}(s \log s + s^{1/2} \epsilon^{-3/2})$ evaluations of $\nabla f_i$ and $\nabla c$, respectively. Moreover, { with a suitable choice of $\{\rho_k\}$,}  it can compute an $\epsilon$-EFSO solution of~\eqref{prob1} using $\widetilde{\mathcal{O}}(s \log s + s^{1/2} \epsilon^{-1})$ evaluations of $\nabla f_i$ and $\nabla c$, respectively.  { These are the best-known first-order oracle complexities in the literature for computing an $\epsilon$-SFSO solution and an $\epsilon$-EFSO solution of problem~\eqref{prob1}, and they are achieved by the SVRFP method. However, there is a significant gap between these two complexity bounds, mainly due to the following two facts established in~\cite{lu2025first}:
\begin{itemize}
\item 
The solution sequence $\{x_k\}$ generated by the SVRFP method satisfies 
\[
\dist\big(c(x_k),-\cK\big)=\cO(\rho_k^{-1/2}), \qquad \bbE[\dist(c(x_k), -\mathcal{K})] = \cO\Big(\rho_k^{-1/2} \sqrt{\bbE[F_{\rho_k}(x_k)-F_{\rho_k}^*]}\,\Big),
\]
 where $F_{\rho_k}(\cdot)$ and $F_{\rho_k}^*$ denote the objective function and the optimal value of 
problem  \eqref{penalty-prob}, respectively. 
\item
The sequence $\{\bbE[F_{\rho_k}(x_k)-F_{\rho_k}^*]\}$ converges to zero, at a rate 
depending on the choice of $\{\rho_k\}$ . 
\end{itemize}
It follows that $\{\dist\big(c(x_k),-\cK\big)\}$ converges to zero at a slower rate than $\{\bbE[\dist(c(x_k), -\mathcal{K})]\}$. 
As a result, the SVRFP method has a worse complexity bound for finding an $\epsilon$-SFSO solution than for finding an $\epsilon$-EFSO solution of problem~\eqref{prob1}.

As discussed above, the SVRFP method is essentially an ``infeasible'' (i.e., exterior) optimization method. 
It drives the iterates to gradually satisfy the constraints at a relatively slow rate.  This behavior leads to the aforementioned significant gap in first-order oracle complexity between finding an $\epsilon$-SFSO solution 
and an $\epsilon$-EFSO solution of~\eqref{prob1}. Motivated by this observation, we aim to  narrow this gap by proposing a ``feasible'' method---an inexact SVRG method---in which all iterates nearly satisfy the constraints.}

Specifically, observe that problem \eqref{prob1} can be reformulated as
\beq \label{prob-reform}
\min_{x} \big\{f(x) +\psi(x)+\delta_{S}(x)\big\},
\eeq 
where $S=\{x\in\rr^n: c(x)\in-\cK\}$, and $\delta_{S}(\cdot)$ denotes the indicator function of the set $S$. 
While problem \eqref{prob-reform} is structurally similar to \eqref{fun-f},  the accelerated (exact) SVRG method in \cite[Algorithm~1]{lan2019unified} is not directly applicable, as its development and analysis require exact evaluation of the proximal operator of $\gamma[\psi(\cdot) + \delta_S(\cdot)]$ for $\gamma > 0$.  However, this proximal operator cannot be  computed exactly in general. Indeed, computing it amounts to solving the constrained optimization problem
\beq \label{prox-prob}
\begin{aligned}
\min_x&\ \, \Big\{\psi(x)+\frac{1}{2\gamma} \|x-u\|^2\Big\} \\
\text{s.t.}&\ \, c(x)\in -\cK
\end{aligned}
\eeq
for any given $u\in\rr^n$. Since the constraint set in~\eqref{prox-prob} is typically non-simple, computing its exact optimal solution is challenging. To address this issue, we propose an inexact SVRG method (Algorithm~\ref{alg2}) for solving~\eqref{prob-reform} (or equivalently, \eqref{prob1}). Our approach follows a similar framework to~\cite[Algorithm~1]{lan2019unified}, but requires only an \emph{inexact} solution to problem~\eqref{prox-prob}---particularly an $\epsilon$-Karush-Kuhn-Tucker ($\epsilon$-KKT) solution $(x,\lambda)$ satisfying
\begin{equation}\label{ekkt}
    \dist\big(0,(x-u)/\gamma+\nabla c(x)\lambda+\partial\psi(x)\big)\leq\epsilon, \qquad \dist\big(c(x),\mathcal{N}_{\cK^*}(\lambda)\big)\leq\epsilon
\end{equation}
for some $\epsilon>0$. Such an $\epsilon$-KKT solution can be computed efficiently using a first-order proximal augmented Lagrangian (AL) method (see Algorithm~\ref{finitesum:alg1}). We show that the proposed inexact SVRG method can find an $\epsilon$-SFSO solution of  \eqref{prob1} using 
\[
\cO\big(s\log s+s^{1/2}\epsilon^{-1/2}\big)\quad  \text{and} \quad  \widetilde\cO\big(s (\log s)^{3/2}\epsilon^{-1/2} + s^{1/4}\epsilon^{-5/4}\big)
\]  
evaluations of $\nabla f_i$ and $\nabla c$, respectively. { These results significantly improve upon the best-known first-order oracle complexity for computing an $\epsilon$-SFSO solution achieved by~\cite[Algorithm~2]{lu2025first}, and substantially reduce the gap between the complexities for obtaining $\epsilon$-SFSO and $\epsilon$-EFSO solutions of~\eqref{prob1}.}

The main contributions of this paper are summarized as follows.

\begin{itemize}
\item We propose an inexact SVRG method (Algorithm~\ref{alg2}) for solving problem~\eqref{prob1}. This method computes an $\epsilon$-SFSO solution---a notion stronger than the commonly studied  $\epsilon$-EFSO solution---using $\cO\big(s\log s + s^{1/2}\epsilon^{-1/2}\big)$ and $\widetilde{\cO}\big(s (\log s)^{3/2}
\epsilon^{-1/2} + s^{1/4}\epsilon^{-5/4}\big)$ evaluations of $\nabla f_i$ and $\nabla c$, respectively. These complexity bounds significantly improve the previously best-known complexity bound $\cO\big(s \log s + s^{1/2} \epsilon^{-3/2}\big)$ established in \cite{lu2025first}. 

\item { In the special case where $s=1$, the proposed method reduces to a novel deterministic first-order method and can compute an $\epsilon$-optimal solution of~\eqref{prob1} using $\cO(\epsilon^{-1/2})$ evaluations of $\nabla f_1$ and $\widetilde{\cO}(\epsilon^{-5/4})$ evaluations of $\nabla c$ and projections onto $-\cK$. This yields a significant improvement in the gradient evaluation complexity with respect to $f_1$ over the best-known first-order methods (e.g., \cite{kovalev2022first,lu2023iteration}), which require $\cO(\epsilon^{-1})$ or $\widetilde{\cO}(\epsilon^{-1})$ evaluations of $\nabla f_1$, $\nabla c$, and projections onto $-\cK$. Consequently, when evaluating $\nabla c$ and projecting onto $-\cK$ are cheaper than evaluating $\nabla f_1$, the proposed method can achieve superior practical performance.}
\end{itemize}

The rest of the paper is organized as follows. Subsection~\ref{sec:notation} introduces the notation and assumptions used throughout the paper. Section~\ref{sec:ALmethod} presents a first-order proximal AL method for solving constrained strongly convex optimization problems, along with its convergence results. In Section~\ref{sec:finite-sum}, we propose an inexact SVRG method for solving problem~\eqref{prob1} and establish its convergence guarantees. Preliminary numerical results are presented in Section~\ref{sec:results} to illustrate the performance of the proposed method. The proofs of the main theoretical results are provided in Section~\ref{sec:proof}. { Finally, concluding remarks are given in Section~\ref{sec:conclude}.}

\subsection{Notation and assumptions}\label{sec:notation}

The following notation will be used throughout this paper. Let $\rr^n$ denote the $n$-dimensional Euclidean space. We use $\langle\cdot, \cdot \rangle$ to denote the standard inner product and $\|\cdot\|$ to denote the Euclidean norm. For any $w \in \rr$, let $\lceil w \rceil$ denote the smallest integer greater than or equal to $w$. The notations $[\cdot]_+$ and $[\cdot]_-$ denote the projection operator onto the nonnegative orthant and nonpositive orthant, respectively.

A mapping $\phi$ is said to be \emph{$L_{\phi}$-Lipschitz continuous} on a set $\Omega$ if $\|\phi(x)-\phi(y)\| \leq L_{\phi} \|x-y\|$ for all $x,y\in \Omega$. It is said to be \emph{$L_{\nabla\phi}$-smooth} on $\Omega$ if $\|\nabla\phi(x)-\nabla\phi(y)\| \leq L_{\nabla\phi} \|x-y\|$ for all $x,y\in \Omega$, where $\nabla \phi$ denotes the transpose of the Jacobian of $\phi$. 

Given a nonempty closed convex set $\Omega \subseteq \rr^n$, let $\inte(\Omega)$ denote its interior, $\dist(x,\Omega)$ the Euclidean distance from $x$ to $\Omega$, and $\Pi_\Omega(x)$ the Euclidean projection of $x$ onto $\Omega$. The indicator function of $\Omega$ is denoted by $\delta_\Omega$, and the normal cone to $\Omega$ at $x \in \Omega$ is denoted by $\cN_\Omega(x)$. For a closed convex cone $\cK \subseteq \rr^n$, its dual cone is denoted by $\cK^*$, that is, $\cK^* = \{y \in \rr^n: \langle y, x\rangle \geq 0 \,\,\, \forall x \in \cK\}$.  Finally, let $\cB$ denote the unit Euclidean ball, i.e., $\cB = \{x: \|x\| \leq 1\}$.

Throughout this paper, we impose the following additional assumptions on problem~\eqref{prob1}.
\begin{assumption} \label{a1}
\begin{enumerate}[label=(\roman*)] 
\item The set $X$ (i.e., the domain of $\psi$) is compact, and $\psi$ is continuous on $X$.
\item Each function $f_i$ is $L_{f_i}$-Lipschitz continuous and $L_{\nabla f_i}$-smooth on $X$ for all $i = 1, 2, \dots, s$.
\item The mapping $c$ is $L_c$-Lipschitz continuous and $L_{\nabla c}$-smooth on $X$, and satisfies $\|c(x)\| \leq C_c$ for all $x \in X$.
\item There exists $\xint \in X$ such that $c(\xint) \in \operatorname{int}(-\cK)$; that is, there exists $D_c > 0$ such that $c(\xint) + D_c \cB \subset -\cK$.
\item The projection onto $\cK^*$ can be exactly computed.
\end{enumerate}
\end{assumption}

In addition, for notational convenience, we define
\begin{align}
& C := L_{\nabla c} C_c + L_c^2, \qquad D_X := \max\Big\{ \max_{x,y \in X} \|x - y\|,\, 1 \Big\}, \label{DX} \\
& D_F := \max_{x,y \in X} \{F(x) - F(y)\}, \qquad D_\psi := \max_{x,y \in X} \{\psi(x) - \psi(y)\}, \label{DF} \\
& f(x) := \frac{1}{s} \sum_{i=1}^s f_i(x), \qquad L_{f} := \frac{1}{s} \sum_{i=1}^s L_{f_i}, \qquad L_{\nabla f} := \frac{1}{s} \sum_{i=1}^s L_{\nabla f_i},  \label{Lf} \\
& \Lambda := \frac{5L_{\nabla f}D_X^2/2+L_f D_X+D_\psi}{D_c}. \label{Lambda} 
\end{align}
It then follows from Assumption~\ref{a1} that $C$, $D_X$, $D_F$, and $D_\psi$ are finite. Moreover, $f$ is $L_{\nabla f}$-smooth on $X$, and
\begin{equation}\label{f-bnd}
    \|\nabla f(x)\| \leq L_{f} \qquad \forall x \in X.
\end{equation}

\section{A first-order proximal augmented Lagrangian method for  constrained strongly convex optimization} \label{sec:ALmethod}

In this section, we present a first-order proximal augmented Lagrangian (FPAL) method for solving a constrained strongly convex  optimization problem, which will serve as a subproblem solver for problem \eqref{prob1}. Specifically, we consider the problem
\beq\label{primal}
\begin{aligned}
\barF^*=\min_{x} &\  \{\barF(x):= \barf (x) +\psi(x)\}\\
\mbox{s.t.}& \ c(x)\in -\cK,
\end{aligned}
\eeq
where $\psi$, $c$, and  $\cK$  are defined in Section \ref{sec:introduction} and satisfy Assumption \ref{a1}. In addition, we impose the following assumption on problem \eqref{primal}. 

\begin{assumption} \label{finite-sum:a2}
The function $\barf$ is $L_{\nabla \barf}$-smooth and $\mu$-strongly convex  on $X$, where $0 < \mu \leq L_{\nabla \barf}$ and $X = \dom\,\psi$ is a compact convex set.
\end{assumption}

Under Assumptions~\ref{a1} and~\ref{finite-sum:a2},  strong duality holds for problem~\eqref{primal} and its Lagrangian dual problem, which is given by
\begin{equation}\label{dual}
\max_{\lambda \in \cK^*} \min_{x} \left\{ \barf(x) + \langle \lambda, c(x) \rangle + \psi(x) \right\}.
\end{equation}
As a result, there exists a pair of optimal solutions $(x^*, \lambda^*) \in X \times \cK^*$ to problems~\eqref{primal} and~\eqref{dual} that satisfy the following Karush-Kuhn-Tucker (KKT) conditions:
\[
0 \in \nabla \barf(x^*)+\nabla c(x^*)\lambda^*+\partial\psi(x^*), \qquad  c(x^*) \in \cN_{\cK^*}(\lambda^*).
\]
Since computing a pair of \emph{exact} optimal solutions $(x^*, \lambda^*)$ is generally intractable, we instead aim to compute an approximate solution, particularly an $\epsilon$-KKT solution as defined below.
  
\begin{defi}
For any $\epsilon>0$,  a pair $(x,\lambda)\in X\times\cK^*$ is called an $\epsilon$-KKT solution of problem \eqref{primal} and its dual if it satisfies
\begin{equation*}
    \dist\big(0,\nabla \barf(x)+\nabla c(x)\lambda+\partial\psi(x)\big)\leq\epsilon, \quad \dist\big(c(x),\cN_{\cK^*}(\lambda)\big)\leq\epsilon.
\end{equation*}
\end{defi}

We next adopt the FPAL method proposed in \cite[Algorithm 2]{lu2023iteration} to compute an $\epsilon$-KKT solution of problem \eqref{primal} and its dual. While the original method was developed for the case where $\barf$ is convex (but not necessarily strongly convex), we apply it to the strongly convex setting and derive improved complexity results by leveraging the strong convexity of $\barf$. 

Specifically, let $(x_k, \lambda_k) \in X \times \cK^*$ denote a pair of approximate solutions generated at the $(k-1)$th iteration. At the $k$th iteration, the FPAL method applies an accelerated proximal gradient (APG) method to approximately solve the proximal augmented Lagrangian subproblem:
\begin{equation}\label{PAL-subprob}
   \min_x\Big\{\overline \cL(x,\lambda_k;\rho_k) :=  \cL(x,\lambda_k;\rho_k)+\frac{1}{2\rho_k} \|x-x_k\|^2\Big\}
\end{equation}
for some penalty parameter $\rho_k > 0$,  where  $ \cL(x,\lambda;\rho)$ is the augmented Lagrangian (AL) function defined as
\[
    \cL(x,\lambda;\rho)=\barf(x)+\frac{1}{2\rho}\Big(\dist^2\big(\lambda+\rho c(x),-\cK\big)-\|\lambda\|^2\Big)+\psi(x).
\]
The Lagrange multiplier $\lambda_{k+1}$ is then updated via the standard dual update scheme (see line 3 of Algorithm~\ref{finitesum:alg1}). 

We now present the FPAL method in Algorithm \ref{finitesum:alg1} for solving problem \eqref{primal}.

\begin{algorithm}[H]
\caption{A first-order proximal augmented Lagrangian method for problem \eqref{primal}}\label{finitesum:alg1}
\begin{algorithmic}[1]
\REQUIRE $\bepsilon > 0$, $(x_0, \lambda_0) \in  X \times \cK^*$, $\rho_0 \geq  1$, $0 < \eta_0 \leq 1$, $\zeta > 1$, $0 < \sigma < 1/\zeta$, $\rho_k= \rho_0 \zeta^k$, $\eta_k = \eta_0 \sigma^k$ for all $k \geq 0$. 
\FOR{$k = 0, 1, \dots$}
    \STATE Apply Algorithm \ref{alg4} in Appendix \ref{appendix} with { $\tf(\cdot) \gets \overline\cL(\cdot,\lambda_k;\rho_k)-\psi(\cdot)$,} $\varepsilon \gets \eta_k$, $\mu_\tf \gets \mu + \rho_k^{-1}$, $\tx_0 \gets x_k$ to the proximal AL subproblem \eqref{PAL-subprob},  and denote its output by $x_{k+1}$.
    \STATE Set $\lambda_{k+1} = \Pi_{\mathcal{K}^*} \left( \lambda_k + \rho_k c(x_{k+1}) \right)$.
    \IF{$\rho_k^{-1} \left\| (x_{k+1}, \lambda_{k+1}) - (x_k, \lambda_k) \right\| +\eta_k \leq \bepsilon$}
        \STATE Terminate the algorithm and output $(x_{k+1}, \lambda_{k+1})$.
    \ENDIF
\ENDFOR
\end{algorithmic}
\end{algorithm}

To analyze the complexity of Algorithm~\ref{finitesum:alg1} for computing an $\bepsilon$-KKT solution of problem~\eqref{primal} and its dual, we introduce some additional notation and state a property of Algorithm~\ref{finitesum:alg1} that will be used in the subsequent analysis.

Let $(x^*,\lambda^*)$ be an arbitrary pair of optimal solutions of \eqref{primal} and \eqref{dual}. Define 
\beq \label{alg1-para}
\begin{aligned}
&  D= \max\big\{\|(x_0,\lambda_0)-(x^*,\lambda^*)\|, 1\big\}, \quad  B = L_{\nabla \barf} + L_{\nabla c}\|\lambda^*\| + L_{\nabla c}D,   \\
&\bar{C} = C \rho_0+B+L_{\nabla c}\rho_0\eta_0/(1-\zeta\sigma)+1, 
\end{aligned}
\eeq
where the constants $C$, $L_{\nabla c}$, and $L_{\nabla \barf}$ are defined in Section \ref{sec:introduction} and Assumption \ref{finite-sum:a2}, respectively.  In addition, it  follows from  Lemmas 3 and 5 in \cite{lu2023iteration} that 
\begin{align}
&\|(x_{k+1},\lambda_{k+1}) - (x_k,\lambda_k)\| \leq \|(x_0,\lambda_0) - (x^*, \lambda^*)\| + \sum_{i=0}^k\rho_i\eta_i   \qquad \forall k \geq 0,  \label{iter-diff} \\
&\|(x_k,\lambda_k) - (x^*,\lambda^*)\| \leq \|(x_0,\lambda_0) - (x^*, \lambda^*)\| + \sum_{i=0}^{k-1}\rho_i\eta_i  \qquad \forall k \geq 0.  \label{iter-bound}
\end{align}

We are now ready to present the complexity result of Algorithm~\ref{finitesum:alg1} for computing an $\bepsilon$-KKT solution to problem~\eqref{primal} and its dual. The proof is deferred to Subsection \ref{proof:thm1}.

\begin{thm}\label{thm:sub-complexity}
 Let $D_X$, $D$, $\bar{C}$, $\bepsilon$, $\rho_0$, $\eta_0$, $\zeta$, $\sigma$, $\mu$ be given in \eqref{DX}, \eqref{alg1-para}, Algorithm \ref{finitesum:alg1}, and Assumption \ref{finite-sum:a2}, respectively,  and let
\begin{align}
 &\bK=\left\lceil\max\left\{1,  \log_\zeta\frac{2(D+\rho_0\eta_0)}{(1-\zeta\sigma)\bepsilon}\right\}\right\rceil,  \label{bK} \\
&\bN =\frac{16\zeta\sqrt{\bar{C}}}{\sqrt{\mu}(\sqrt{\zeta}-1)}\log\frac{2\zeta  \bar{C} D_X}{\eta_0\sigma}\left\lceil\max\left\{1,  \log_\zeta\frac{2(D+\rho_0\eta_0)}{(1-\zeta\sigma)\bepsilon}\right\}\right\rceil\max\left\{1,  \sqrt{\frac{2(D+\rho_0\eta_0)}{(1-\zeta\sigma)\bepsilon}}\right\}. \label{bN}
\end{align}
Then,  Algorithm \ref{finitesum:alg1} terminates and outputs an $\bepsilon$-KKT solution to problem~\eqref{primal} and its dual in at most $\bK+1$ iterations. Moreover, the number of evaluations of $\nabla \barf$ and $\nabla c$ performed by Algorithm \ref{finitesum:alg1} is at most $\bN$, respectively.
\end{thm}

\section{An inexact stochastic variance-reduced  gradient method for problem \eqref{prob1}}\label{sec:finite-sum}

In this section, we propose an inexact stochastic variance-reduced  gradient (SVRG) method for solving problem \eqref{prob1}.

Observe that when $c=0$, problem \eqref{prob1} reduces to the unconstrained finite-sum optimization problem  \eqref{fun-f}, which has been extensively studied in the literature (see, e.g.,  \cite{allen2018katyusha,defazio2014saga,hendrikx2021optimal,lan2019unified,
schmidt2017minimizing,shalev2013stochastic,song2020variance,xiao2014proximal,zhou2019lower}).
Notably, an accelerated SVRG method was proposed in \cite[Algorithm 1]{lan2019unified} for solving~\eqref{fun-f}. This method can be viewed as a restarted accelerated stochastic gradient scheme with  two nested loops: the outer loop updates a reference point $\tilde x$ and the initialization points $x_0$, $z_0$ for the inner loop, while the inner loop performs a sequence of accelerated stochastic gradient steps. At each inner iteration, an SVRG estimator of $f$ is computed using the full gradient $\nabla f(\tilde x)$ at the reference point $\tilde x$, along with the sampled gradients $\nabla f_i(\tilde x)$ and $\nabla f_i(y)$, where $y$ is updated at each inner iteration, and  $i\in \{1, 2, \dots, s\}$ is drawn according to a prescribed probability distribution.

Note that problem \eqref{prob1} can be reformulated as \eqref{prob-reform}. While problem \eqref{prob-reform} shares a similar structure with \eqref{fun-f},  the accelerated SVRG method in \cite[Algorithm 1]{lan2019unified} is not  directly applicable to \eqref{prob-reform} because the proximal operator of $\gamma [\psi(\cdot)+\delta_{S}(\cdot)]$ is not exactly computable for general $\gamma>0$. In fact, evaluating this proximal operator requires solving problem \eqref{prox-prob}, which typically involves a non-simple constraint set. As a result, computing an exact solution to~\eqref{prox-prob} is generally challenging. 
To address this issue, we propose an inexact SVRG method for solving problem \eqref{prob-reform} (equivalently, problem \eqref{prob1}), which follows a similar framework to \cite[Algorithm 1]{lan2019unified} but only requires  an inexact solution to problem \eqref{prox-prob}---particularly an $\epsilon$-KKT solution $(x,\lambda)$ satisfying \eqref{ekkt} for some suitable $\epsilon>0$. As discussed in Section \ref{sec:ALmethod}, such an $\epsilon$-KKT solution can be efficiently computed by Algorithm \ref{finitesum:alg1}. 

We now present an inexact SVRG method for solving problem \eqref{prob1} in Algorithm \ref{alg2}. 

\begin{algorithm}[H]
\caption{An inexact SVRG method for problem \eqref{prob1}}\label{alg2}
\begin{algorithmic}[1]
\REQUIRE $\epsilon>0$, $R\geq 0$, $\{\alpha_k\}\subset(0, 1]$, $\{p_k\}\subset(0,1]$, $\{\gamma_k\}\subset (0,\infty)$, $\{T_k\}$, $\{\theta_t\}$, a probability distribution $\mathcal{Q}=\{q_1, q_2,\dots, q_s\}$.
\STATE Call Algorithm \ref{alg3} in Appendix \ref{appendix} to solve the problem $\min_{x \in X} \dist^2\big(c(x),-\cK\big)$ to find $x^0\in X$ such that $\dist\big(c(x^0),-\cK\big)\leq\epsilon$.
\STATE Set $\tilde x_1=x^0$, $\tilde z_1=x^0$, and $\tilde\lambda_1=0$.
\FOR{$k=1,2,\dots$}
\STATE Set $x_0=\tilde x_{k}$, $z_0=\tilde z_{k}$,  $\lambda_0=\tilde\lambda_k$, and $\tilde g_{k} = \nabla f(\tilde x_{k})$. 
\FOR{$t=1, 2, \dots, T_k$}
\STATE Sample $i_t\in \{1, 2, \dots, s\}$ according to $\mathcal{Q}$.
\STATE $y_t=(1-\alpha_k-p_k)x_{t-1}+\alpha_kz_{t-1}+p_k\tilde x_{k}$.
\STATE $g_t=(\nabla f_{i_t}(y_t)-\nabla f_{i_t}(\tilde x_{k}))/(q_{i_t}s)+\tilde g_{k}$.
\STATE Call Algorithm \ref{finitesum:alg1} with $\bepsilon \gets \epsilon$, $\mu \gets 1/\gamma_k$, $(x_0,\lambda_0) \gets (z_{t-1},\lambda_{t-1} \min\{R/\|\lambda_{t-1}\|, 1\})$ to solve the subproblem
\begin{equation}\label{subproblem}
\begin{aligned}
    \min_x  & \ \left\{F_{k,t}(x) :=\frac{1}{2\gamma_k} \|x-z_{t-1}+\gamma_kg_t\|^2+\psi(x)\right\}  \\
 \mbox{s.t.} & \ c(x) \in -\cK,
\end{aligned}
\end{equation}
and denote its output as $(z_t,\lambda_t)$.\footnotemark
\STATE $x_t=(1-\alpha_k-p_k)x_{t-1}+\alpha_kz_t+p_k\tilde x_k$.
\ENDFOR
\STATE Set $\tilde z_{k+1}=z_{T_k}$, $\tilde x_{k+1}=\sum_{t=1}^{T_k}(\theta_tx_t)/\sum_{t=1}^{T_k}\theta_t$, and $\tilde\lambda_{k+1}=\lambda_{T_k}$.
\ENDFOR
\end{algorithmic}
\end{algorithm}
 \footnotetext{By convention, we set $0/0=0$ and $1/0=\infty$.}

\begin{rem} 
(i) Algorithm~\ref{alg2} shares some similarities with the gradient sliding (GS) method~\cite{lan2016gradient}. In particular, the latter is a double-loop first-order method for solving problems of the form 
\beq \label{gradslide-prob}
\min_x \{g(x)+h(x)+\psi(x)\},
\eeq
 where $g$ is a Lipschitz-smooth convex function, $\psi$ is the same as in~\eqref{prob1}, and $h$ is a nonsmooth convex function such that $\operatorname{ri}(\dom\, h) \supseteq \dom\,\psi$.\footnote{Here, $\operatorname{ri}(\dom\, h)$ denotes the relative interior of $\dom\,h$, and the latter assumption ensures that $\partial h(x) \neq \emptyset$ for all $x \in \dom\,\psi$, which is required for applying a proximal subgradient method to solve subproblems of the form~\eqref{gradslide-subprob}.} The outer loop of GS performs an inexact accelerated proximal gradient method, whose main computational effort lies in solving subproblems of the form
\beq \label{gradslide-subprob}
\min_x \Big\{\frac{1}{2\gamma} \|x-u\|^2 + h(x) + \psi(x)\Big\}
\eeq
for $\gamma>0$ and a vector $u$ obtained using the gradient of $g$ at some point. The inner loop of GS applies a proximal subgradient method to solve these subproblems, where the main computational cost lies in evaluating the proximal operator of $\psi$ at points obtained using subgradients of $h$.

However, although problem~\eqref{prob-reform} shares the same form as~\eqref{gradslide-prob}, it does not belong to this problem class, since $\operatorname{ri}(\dom\,\delta_S) = \operatorname{ri}(S)$ typically does not contain $\dom\,\psi$. Consequently, the GS method~\cite{lan2016gradient} is not applicable to problem~\eqref{prob-reform} (or equivalently,~\eqref{prob1}).

(ii) The parameter $R$ ensures the boundedness of $\lambda_t$, where $\lambda_t$ is an approximate optimal  Lagrange multiplier of subproblem~\eqref{subproblem} computed in Line~9 of Algorithm~\ref{alg2}. Indeed, as will be shown in Lemmas~\ref{l-subbound} and~\ref{iter-property}(iii), it holds that
\beq \label{lambda-bnd}
\|\lambda_t\| \leq  D_X + R+2\Lambda+\rho_0\eta_0/(1-\zeta\sigma),
\eeq 
where $D_X$ and $\Lambda$ are defined in~\eqref{DX} and~\eqref{Lambda}, respectively, and $\rho_0$, $\zeta$, and $\sigma$ are input parameters of Algorithm~\ref{finitesum:alg1}. 

In addition, note from Line~9 of Algorithm~\ref{alg2} that $\lambda_0 = \alpha \lambda_{t-1}$, where $\alpha = \min\{R/\|\lambda_{t-1}\|, 1\}$, and $\lambda_0$ serves as the initial guess for the optimal Lagrange multiplier of~\eqref{subproblem}. Clearly, if $R=0$, then $\lambda_0=0$, which corresponds to a static initialization for Algorithm \ref{finitesum:alg1}. On the other hand, by~\eqref{lambda-bnd} and the definition of $\alpha$, we have
\[
\frac{R}{D_X + R+2\Lambda+\rho_0\eta_0/(1-\zeta\sigma)} \leq \alpha \leq 1.
\]
Consequently, when $R$ is large, $\alpha \approx 1$ and hence $\lambda_0 \approx \lambda_{t-1}$, which provides a dynamic (warm-start) initialization for Algorithm \ref{finitesum:alg1}.

As will be shown in Theorem~\ref{thm:finite-sum}, setting $R=0$ (static initialization) yields the lowest first-order oracle complexity bounds among all possible choices of $R$. However, numerical experiments indicate that larger values of $R$ (warm-start) often lead to faster convergence in practice.
\end{rem}

The following theorem presents the convergence results and provides an estimate of the computational cost of Algorithm~\ref{alg2} in terms of the number of evaluations of $\nabla f_i$ and $\nabla c$. The proof is deferred to Subsection~\ref{proof:thm2}.

\begin{thm}\label{thm:finite-sum}
 Suppose that Assumption \ref{a1} holds. Let $D_X$, $D_F$, $L_{\nabla f}$, $\Lambda$, $C_c$,    $L_c$, $L_{\nabla c}$, $L_{\nabla f_i}$, $\epsilon$,   $\rho_0$, $\eta_0$, $\zeta$, $\sigma$, and $R$ be  given  in \eqref{DX}, \eqref{DF}, \eqref{Lf}, \eqref{Lambda}, Assumption \ref{a1}, and Algorithms \ref{finitesum:alg1} and \ref{alg2}, respectively, and let 
\begin{align}
C_1 &= 2 D_X^{3/2} \sqrt{L_c^2+L_{\nabla c}C_c}, \nn \\ 
C_2 &=  (L_{\nabla c}C_c+L_c^2)\rho_0+3L_{\nabla f}/2  + L_{\nabla c}(D_X + R +2 \Lambda)+L_{\nabla c}\rho_0\eta_0/(1-\zeta\sigma)+1,\label{C2} \\
C_3 &= L_{\nabla f}^{-\frac{1}{2}} \frac{16\zeta\sqrt{C_2}}{\sqrt{\zeta}-1} \log\frac{2\zeta  C_2 D_X}{\eta_0\sigma}, \qquad C_4 = \frac{2(D_X + R + \Lambda+\rho_0\eta_0)}{1-\zeta\sigma},  \label{C3} \\
C_5& =3D_X + 2R+4\Lambda+2\rho_0\eta_0/(1-\zeta\sigma), \qquad C_6=2D_F+3L_{\nabla f}D_X^2/2. \label{c45}
\end{align}
Assume that $\{\tx_k\}$ is generated by Algorithm \ref{alg2} with the parameters chosen as
\begin{equation}\label{def-para}
\begin{aligned}
    & k_0=\lfloor \log_2 s \rfloor+1, \quad  \alpha_k=\left\{
    \begin{array}{ll}
      \frac{1}{2}   & \text{if} \  \ 1\leq k\leq k_0, \\
       \frac{2}{k-k_0+4}  & \text{if} \  \ k>k_0,
    \end{array}\right.  \quad  p_k=\frac{1}{2},  \quad \gamma_k=\frac{1}{3L_{\nabla f}\alpha_k},\footnotemark \\
&       q_i=\frac{L_{\nabla f_i}}{\sum_{j=1}^s L_{\nabla f_j}}, \quad T_k=\left\{\begin{array}{ll}
       2^{k-1}  & \text{if} \  \ 1\leq k\leq k_0, \\
       2^{k_0-1}  & \text{if} \  \ k>k_0,
       \end{array}\right. 
    \quad \theta_t=\left\{
    \begin{array}{ll}
      \frac{\gamma_k}{\alpha_k}(\alpha_k+p_k)   & \text{if} \  \ 1\leq t\leq T_k-1, \\
       \frac{\gamma_k}{\alpha_k}  & \text{if} \  \ t=T_k.
    \end{array}\right.  
\end{aligned}  
\end{equation}
\footnotetext{{ We choose $p_k = 1/2$ for simplicity in the analysis of Algorithm~\ref{alg2}. However, dynamic choices are also possible, provided they satisfy the conditions specified in Lemma~\ref{l-iter2}.}}
Let $N_f$ and $N_c$ denote the number of evaluations of $\nabla f_i$ and $\nabla c$  during the $k$ outer iterations of Algorithm \ref{alg2} for all $k \geq 1$, respectively.
Then,  for all $k>k_0$,  it holds that
\begin{align}
& \dist\big(c(\tx_k),-\cK\big) \leq \epsilon, \qquad F(\tilde x_k)-F^* \geq -D_\lambda \epsilon, \qquad  \bbE[F(\tilde x_k)-F^*]\leq
      \frac{16C_6}{s(k-k_0+3)^2}+8C_5\epsilon,  \label{xk-prop} \\
& N_f \leq2sk,  \qquad N_c\leq C_1\epsilon^{-1} +2 s C_3  \left\lceil \max\left\{1, \log_\zeta\frac{C_4}{\epsilon}\right\}\right\rceil \max\left\{1,  \sqrt{\frac{C_4}{\epsilon}} \right\}k^{\frac{3}{2}},  \label{N-bound}
\end{align}
where $D_\lambda$ denotes the distance between the origin and the set of optimal Lagrange multipliers of \eqref{prob1}.
\end{thm}

\begin{rem}\label{rem2}
\begin{itemize}
\item[(i)]
Letting $k=\ceil {k_0+s^{-1/2}\epsilon^{-1/2}}$, one can observe from \eqref{xk-prop} that 
\begin{align*}
& \dist(c(\tx_k),-\cK) \leq \epsilon, \qquad F(\tilde x_k)-F^* \geq -D_\lambda \epsilon, \qquad  \bbE[F(\tilde x_k)-F^*]\leq 8(C_5+2C_6)\epsilon, \\
& N_f=\cO\big(s\log s+s^{1/2}\epsilon^{-1/2}\big), \qquad N_c= \widetilde\cO\big(s (\log s)^{3/2}\epsilon^{-1/2} + s^{1/4}\epsilon^{-5/4}\big).
\end{align*}
As a result, Algorithm \ref{alg2} can find an $\epsilon$-SFSO solution of problem \eqref{prob1} using $\cO\big(s\log s+s^{1/2}\epsilon^{-1/2}\big)$ and $\widetilde\cO\big(s (\log s)^{3/2}\epsilon^{-1/2} + s^{1/4}\epsilon^{-5/4}\big)$ evaluations of $\nabla f_i$ and $\nabla c$, respectively. In contrast, when $\cK$ is the nonnegative orthant, the stochastic variance-reduced first-order penalty method proposed in \cite[Algorithm 2]{lu2025first} requires $\cO(s \log s + s^{1/2} \epsilon^{-3/2})$ evaluations of both $\nabla f_i$ and $\nabla c$ to compute an $\epsilon$-SFSO solution of \eqref{prob1}. Therefore, Algorithm~\ref{alg2} significantly improves the first-order oracle complexity achieved by \cite[Algorithm 2]{lu2025first}.
\item[(ii)] We now consider the special case of problem~\eqref{prob1} with $s=1$. In this case, problem~\eqref{prob1} reduces to a standard convex conic optimization problem, and Algorithm~\ref{alg2} becomes a novel deterministic first-order method. As discussed above, Algorithm~\ref{alg2} can compute an $\epsilon$-optimal solution of~\eqref{prob1} using $\cO(\epsilon^{-1/2})$ and $\widetilde{\cO}(\epsilon^{-5/4})$ evaluations of $\nabla f_1$ and $\nabla c$, respectively. In contrast, the best-known first-order methods---such as the extra anchored gradient method \cite[Algorithm 3]{kovalev2022first} and the first-order proximal augmented Lagrangian method \cite[Algorithm 2]{lu2023iteration}---require $\cO(\epsilon^{-1})$ or $\widetilde{\cO}(\epsilon^{-1})$ evaluations of $\nabla f_1$, $\nabla c$, and projection onto $-\cK$ to obtain an $\epsilon$-optimal solution of~\eqref{prob1}. Compared to these methods, Algorithm~\ref{alg2} significantly reduces the gradient evaluation complexity with respect to $f_1$, at the cost of a higher complexity with respect to $c$. { Note that the number of projections onto $-\cK$ is on the same order as the number of evaluations of $\nabla c$ in Algorithm~\ref{alg2}. Consequently, when computing $\nabla c$ and projecting onto $-\cK$ are cheaper than evaluating $\nabla f_1$, Algorithm~\ref{alg2} can achieve superior practical performance.}
\end{itemize}
\end{rem}

\section{Numerical results} \label{sec:results}
In this section, we conduct preliminary experiments to test the performance of our proposed method (Algorithm \ref{alg2}), and compare it with the stochastic variance-reduced  first-order penalty (SVRFP) method \cite[Algorithm 2]{lu2025first} and the extra anchored gradient (EAG) method \cite[Algorithm 3]{kovalev2022first}. EAG is an optimal first-order method designed for monotone inclusion problems. Since problem~\eqref{prob1} can be equivalently reformulated as a monotone inclusion problem, EAG is well-suited for solving it. The experiments are coded in MATLAB and all the computations are conducted on a desktop with a 2 GHz Intel i9-14900F 24-core processor and 32 GB of RAM. 

\subsection{Fairness-constrained logistic regression}
\label{sec:exp1}

We consider the dataset \texttt{a1a} from LIBSVM \cite{chang2011libsvm}, which consists of approximately $67\%$ male samples and $33\%$ female samples. Due to this imbalance, a classical unconstrained logistic regression model trained on the full dataset typically performs well on the male subgroup but poorly on the female subgroup, indicating potential unfairness. Therefore, incorporating a fairness constraint \cite{diana2021minimax,zafar2017fairness} is necessary. 

To address this issue, we minimize the empirical logistic loss subject to a group fairness constraint, defined as a performance guarantee on a validation set $V$ randomly drawn from the female subgroup of the dataset. The resulting fairness-constrained logistic regression model  is given by
\begin{equation} \label{fair-logistic}
  \min_{w,b} \;
    \frac{1}{s}\sum_{i=1}^{s} \ell(w,b;\, x_i, y_i) + \lambda\|w\|_1
  \quad \mbox{s.t.} \quad
  \frac{1}{|V|}\sum_{j \in V} \ell(w,b;\, x_j, y_j) \leq \delta
\end{equation}
for some $\lambda, \delta>0$, where $x_i \in \rr^n$ and $y_i \in \{-1,1\}$, $i=1,\ldots, s$, are the sample data from \texttt{a1a}, and $\ell$ is the logistic loss function given by $\ell(w,b; x, y) = \log\big(1 + e^{-y (w^Tx+b)}\big)$.

We set $\lambda = 0.01$, $|V| = 100$, and $\delta = (\widehat L_V+L_V^*)/2$ for \eqref{fair-logistic}, where
$\widehat L_V = |V|^{-1}\sum_{j\in V}\ell(\hat w,\hat b;\,x_j,y_j)$ is the validation loss evaluated at $(\hat w,\hat b)$, an optimal solution of the unconstrained counterpart of \eqref{fair-logistic} (i.e., problem \eqref{fair-logistic} with $\delta = \infty$), and 
\begin{equation} \label{LV-star}
L_V^* = \min_{w,b}|V|^{-1}\sum_{j\in V}\ell(w,b,x_j,y_j).
\end{equation}
Clearly, this choice of $\delta$ ensures that problem \eqref{fair-logistic} is feasible. Moreover, since $L_V^*<\widehat L_V$ typically holds, we have $\delta < \widehat L_V$, implying that $(\hat w,\hat b)$ is infeasible for \eqref{fair-logistic}. Consequently, the constraint in \eqref{fair-logistic} is generally active at any optimal solution.

 \begin{figure}[h]
    \centering
    \includegraphics[width=0.8\linewidth]{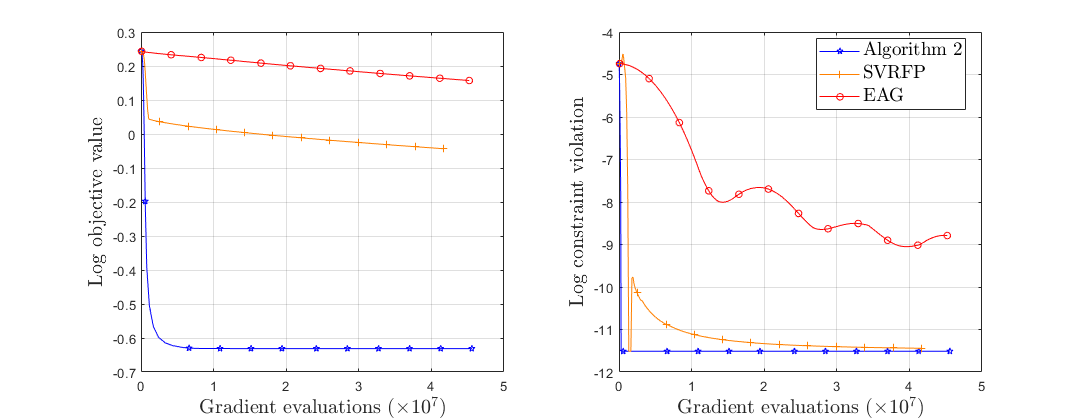}
    \caption{Convergence behavior of objective value and constraint violation for fairness-constrained logistic regression.}
    \label{fig:1}
\end{figure}

To apply Algorithm~\ref{alg2}, SVRFP~\cite{lu2025first}, and EAG~\cite{kovalev2022first} to solve \eqref{fair-logistic}, we specify an initial point and the corresponding parameters. Specifically, we solve problem \eqref{LV-star} to obtain an $10^{-2}$-optimal solution $(w_0, b_0)$, and use it as the initial point for all three methods. 
  The parameters for these methods are chosen as follows:
\begin{itemize}
    \item For Algorithm~\ref{alg2}, we set $(\epsilon,\rho_0,\eta_0,\zeta,\sigma, R) = (0.01, 10, 0.1, 3, 0.3,10)$, and choose $(\alpha_k, p_k, \gamma_k, q_i, T_k, \theta_t)$ as specified in Theorem~\ref{thm:finite-sum}.
    \item For SVRFP, the parameters $(\alpha_k, p_k, \gamma_k, q_i, T_k, \theta_t, \rho_k)$ are selected as in \cite[Theorem 4(i)]{lu2025first}.
    \item For EAG, we select $(\beta_t, \lambda)$ according to \cite[Theorem 2]{kovalev2022first}.
\end{itemize}

We now apply these three methods to solve \eqref{fair-logistic} and compare their performance. In particular,  we evaluate the objective value and constraint violation as the number of gradient evaluations increases for all three methods. The results are presented in Figure~\ref{fig:1}, which shows that the stochastic methods---Algorithm~\ref{alg2} and SVRFP---converge significantly faster than the deterministic method EAG, with Algorithm~\ref{alg2} converging faster than SVRFP in terms of objective value. In addition, Algorithm~\ref{alg2} maintains a consistently low level of constraint violation throughout the iterations, as it is designed to nearly achieve feasibility at each iteration. In contrast, both SVRFP and EAG exhibit an initial increase in constraint violation followed by a gradual decrease.

\subsection{Quadratically constrained least squares regression}
\label{sec:exp2}

We consider the following quadratically constrained least-squares problem \cite{lawson1995solving}:
\begin{equation} \label{q-regress}
  \min_{x} \;
    \frac{1}{s}\sum_{i=1}^{s} \|A_i x - b_i\|^2 + \lambda\|x\|^2
  \quad \mbox{s.t.} \quad
  \|Cx - d\|^2 \leq r^2
\end{equation}
for some $\lambda, r>0$, where $A_i \in \mathbb{R}^{l \times n}$, $b_i \in \mathbb{R}^l$,
$C \in \mathbb{R}^{p \times n}$, and $d \in \mathbb{R}^p$.
Quadratic constraints of this form often arise in signal processing (e.g., power budgets)
and control (e.g., ellipsoidal uncertainty sets), where feasibility is critical.

We generate data for \eqref{q-regress} with $n = 100$, $s = 10{,}000$, $l = 50$, and $p = 20$. Specifically, the data matrices $A_i$ and $C$ have entries drawn independently from
$\mathcal{N}(0, l^{-1})$ and $\mathcal{N}(0, p^{-1})$, respectively. We then generate a sparse vector $\hat x \in \mathbb{R}^n$ with $10\%$ nonzero entries, whose nonzero components are drawn from the standard normal distribution $\mathcal{N}(0,1)$, and set $b_i = A_i \hat x + \xi_i$, where $\xi_i \sim \mathcal{N}(0, 10^{-4}I)$.

We set $\lambda = 0.01$ and $r = 0.5\,\|C\hat x^* - d\|$ for \eqref{q-regress}, where
${\hat x}^* = \big(\sum_{i=1}^s A_i^T A_i + s\lambda I\big)^{-1} \sum_{i=1}^s A_i^T b_i$
is an optimal solution of the unconstrained counterpart of \eqref{q-regress} (i.e., problem \eqref{q-regress} with $r = \infty$), and $d \sim \mathcal{N}(0, I)$ is generated such that $C\hat x^* \neq d$. This condition, together with the choice of $r$, excludes $\hat x^*$ from being an optimal solution of \eqref{q-regress}. Consequently, the constraint in \eqref{q-regress} must be active at any optimal solution.

To apply Algorithm~\ref{alg2}, SVRFP~\cite{lu2025first}, and EAG~\cite{kovalev2022first} to solve \eqref{q-regress}, we specify an initial point and the corresponding parameters. Specifically, we solve  $\min_{x} \|Cx-d\|^2$ to obtain a $10^{-2}$-feasible solution $x_0$ of \eqref{q-regress}, i.e., $\|Cx_0-d\|^2 \leq r^2 +10^{-2}$. We then use $x_0$ as the initial point for all three methods. The parameters for these methods are chosen to be the same as those specified in Subsection \ref{sec:exp1}.

 \begin{figure}[h]
    \centering
    \includegraphics[width=0.8\linewidth]{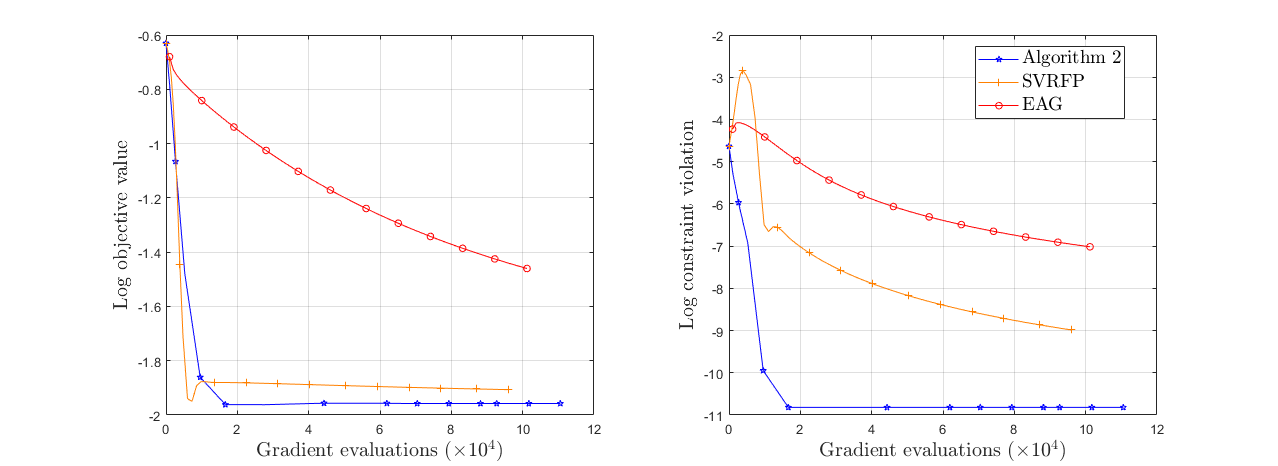}
    \caption{Convergence behavior of objective value and constraint violation for quadratically constrained least squares regression.}
    \label{fig:2}
\end{figure}

We now apply these three methods to solve \eqref{q-regress} and compare their performance. In particular,  we evaluate the objective value and constraint violation as the number of gradient evaluations increases for all three methods. The results are presented in Figure~\ref{fig:2}, which shows that the stochastic methods---Algorithm~\ref{alg2} and SVRFP---converge significantly faster than the deterministic method EAG, with Algorithm~\ref{alg2} converging faster than SVRFP in terms of objective value. In addition, Algorithm~\ref{alg2} achieves a rapid reduction in constraint violation and maintains it at a low level, whereas SVRFP exhibits an initial increase in constraint violation followed by a gradual decrease, and EAG maintains a substantially higher level of violation throughout the process.

\subsection{Linearly constrained LASSO regression}
\label{sec:exp3}

In this subsection, we consider a special case of problem~\eqref{prob1} with $s=1$. In this case, as discussed in Remark~\ref{rem2}(ii), Algorithm~\ref{alg2} reduces to a deterministic first-order method. Moreover, when computing $\nabla c$ and projecting onto $-\cK$ are cheaper than evaluating $\nabla f_1$, Algorithm~\ref{alg2} can achieve superior practical performance compared to SVRFP~\cite{lu2025first} and EAG~\cite{kovalev2022first}. We now provide a numerical example to illustrate this observation.

Specifically, we consider the following linearly constrained LASSO problem \cite{tibshirani1996regression}:
\begin{equation} \label{lasso}
  \min_{x } \;
    \frac{1}{2}\|Ax - b\|^2 + \lambda\|x\|_1
  \quad \mbox{s.t.} \quad
  Cx \leq d
\end{equation}
for some $\lambda>0$, where $A \in \mathbb{R}^{m \times n}$, $b\in \mathbb{R}^m$, $C \in \mathbb{R}^{p \times n}$, and $d\in \rr^p$.  Note that \eqref{lasso} is a special case of problem~\eqref{prob1} with $s=1$, where $f_1(x) = \|Ax - b\|^2/2$, $\psi(x) = \|x\|_1$, $c(x) = Cx - d$, and $\cK = \rr^p_+$. Moreover, computing $\nabla c$ and projecting onto $-\cK$ are typically cheaper than evaluating $\nabla f_1$.

We set $\lambda=0.01$ and generate data for \eqref{lasso} with $m = 2{,}000$, $n = 2{,}000$, and $p = 20$. Specifically, we generate $A$, $b$, and $C$ in the same manner as in Subsection \ref{sec:exp2} such that $C\hx^* \neq 0$, where $\hx^*$ is an optimal solution of the unconstrained counterpart of \eqref{lasso}, i.e., 
$\hx^* \in \argmin_x \{\|Ax-b\|^2/2+\lambda\|x\|_1\}$. We then set $d = 0.5[C\hx^*]_++2[C\hx^*]_-$. This choice of $d$, together with $C\hx^* \neq 0$, ensures that $\hx^*$ is infeasible for \eqref{lasso}. Consequently, at least one of the constraints in \eqref{lasso} must be active at any optimal solution.

To apply Algorithm~\ref{alg2}, SVRFP~\cite{lu2025first}, and EAG~\cite{kovalev2022first} to solve \eqref{lasso}, we specify an initial point and the corresponding parameters. Specifically, we solve the feasibility problem $\min_{x}\|[Cx-d]_+\|^2$ to obtain a $10^{-2}$-feasible solution $x_0$ of \eqref{lasso}, i.e., $\|[Cx_0-d]_+\| \leq 10^{-2}$. We then use $x_0$ as the initial point for all three methods. The parameters for Algorithm~\ref{alg2} and EAG are chosen to be the same as those specified in Subsection \ref{sec:exp1}, while the parameters $(\alpha_k, p_k, \gamma_k, q_i, T_k, \theta_t, \rho_k)$ for SVRFP are selected as in \cite[Theorem 4(ii)]{lu2025first} since it reduces to a deterministic method when $s=1$.

\begin{figure}[h]
    \centering
    \includegraphics[width=0.8\linewidth]{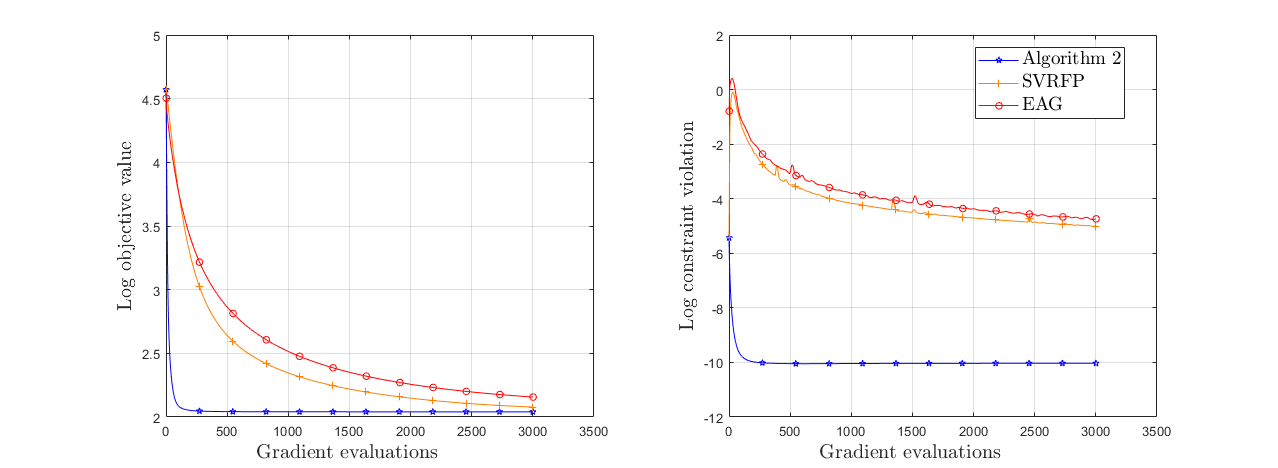}
    \caption{Convergence behavior of objective value and constraint violation for linearly constrained LASSO regression.}
    \label{fig:3}
\end{figure}

We now apply these three methods to solve \eqref{lasso} and compare their performance. Note that all methods reduce to deterministic ones since $s=1$. In particular, we evaluate the objective value and constraint violation as the number of gradient evaluations increases for all three methods. The results are presented in Figure~\ref{fig:3}, which shows that Algorithm~\ref{alg2} converges significantly faster than SVRFP and EAG, while SVRFP and EAG exhibit similar convergence behavior in terms of objective value.  In addition, Algorithm~\ref{alg2} achieves a rapid reduction in constraint violation and maintains it at a low level, whereas SVRFP and EAG exhibit similar but much slower convergence in constraint violation. This phenomenon is consistent with the observation in Remark~\ref{rem2}(ii). Furthermore, the similar performance of SVRFP and EAG is expected, as they have comparable theoretical convergence rates.

\section{Proof of the main results}\label{sec:proof}
In this section, we present a proof of the main results stated in Sections \ref{sec:ALmethod} and \ref{sec:finite-sum}, which are particularly Theorems \ref{thm:sub-complexity} and \ref{thm:finite-sum}. 

\subsection{Proof of the main results in Section \ref{sec:ALmethod}} \label{proof:thm1}

In this subsection, we present a proof of Theorem \ref{thm:sub-complexity}. To this end, we begin by stating a lemma concerning the Lipschitz smoothness of the function $\overline{\mathcal{L}}(\cdot, \lambda_k; \rho_k)$, which was established in \cite[Lemma 7]{lu2023iteration}.

\begin{lemma}[\textbf{\cite[Lemma 7]{lu2023iteration}}] \label{lem:Lk}
Let $\{\rho_k\}$, $\{\eta_k\}$, and $\{\lambda_k\}$ be generated in Algorithm \ref{finitesum:alg1}, and $B$, $C$,  and $L_{\nabla c}$ be given in \eqref{DX}, \eqref{alg1-para}, and Assumption \ref{a1}, respectively. Suppose that Assumptions \ref{a1} and \ref{finite-sum:a2} hold. Then the function $\overline \cL(\cdot,\lambda_k;\rho_k)$ is $L_k$-smooth on $X$ for each $k \geq 0$, where 
$\overline \cL(\cdot,\lambda_k;\rho_k)$ is defined in \eqref{PAL-subprob}, and 
\beq \label{L_k}
     L_k= C \rho_k+B+L_{\nabla c}\sum^{k-1}_{i=0}\rho_i\eta_i+\rho_k^{-1}.
\eeq
\end{lemma}

We are now ready to prove Theorem \ref{thm:sub-complexity}.

\begin{proof}[\textbf{Proof of Theorem~\ref{thm:sub-complexity}}]
We first show that Algorithm \ref{finitesum:alg1} terminates in at most $\bK+1$ outer iterations, where $\bK$ is defined in \eqref{bK}. Indeed, suppose for contradiction that it runs for more than $\bK+1$ iterations. It then follows from Algorithm \ref{finitesum:alg1} that 
\beq \label{non-termination}
\rho_\bK^{-1} \big\| (x_{\bK+1}, \lambda_{\bK+1}) - (x_\bK, \lambda_\bK) \big\| +\eta_\bK > \bepsilon.
\eeq
Notice from Algorithm \ref{finitesum:alg1} that $\zeta>1$ and $\rho_\bK=\rho_0\zeta^\bK$. By these facts and \eqref{bK}, one has 
\beq \label{rho-bound}
    \rho_\bK \geq \rho_0\zeta^\bK\overset{\eqref{bK}}{\geq} \frac{2\rho_0(D+\rho_0\eta_0)}{(1-\zeta\sigma)\bepsilon}.
\eeq
Using this, \eqref{alg1-para}, \eqref{iter-diff},  and the fact that $\rho_0\geq1$, $0<\zeta\sigma<1$, $\rho_k= \rho_0 \zeta^k$, and $\eta_k = \eta_0 \sigma^k$,  we obtain that
\begin{align}
\frac{\| (x_{\bK+1}, \lambda_{\bK+1}) - (x_\bK, \lambda_\bK)\|}{\rho_\bK} & \overset{\eqref{alg1-para}\eqref{iter-diff}}{\leq} \frac{D + \sum_{k=0}^\bK\rho_k\eta_k}{\rho_\bK} \overset{\eqref{rho-bound}}{\leq} \frac{D+\rho_0\eta_0\sum_{k=0}^\bK(\zeta\sigma)^k}{{2\rho_0(D+\rho_0\eta_0)}/{((1-\zeta\sigma)\bepsilon)}}  \nn \\
& \leq  \frac{D+\rho_0\eta_0/(1-\zeta\sigma)}{{2\rho_0(D+\rho_0\eta_0)}/{((1-\zeta\sigma)\bepsilon)}}  =\frac{\bepsilon}{2}\cdot\frac{D(1-\zeta\sigma)+\rho_0\eta_0}{\rho_0(D+\rho_0\eta_0)}\leq\frac{\bepsilon}{2}. \label{relative-change}
\end{align}
In addition, notice that $1<\zeta<\sigma^{-1}$ and $D\geq 1$. By these facts and \eqref{bK}, one has
\begin{equation*}
 \bK\geq\log_\zeta\frac{2(D+\rho_0\eta_0)}{(1-\zeta\sigma)\bepsilon}\geq\log_{\zeta}(2\bepsilon^{-1})\geq\log_{\sigma^{-1}}(2\bepsilon^{-1}) = \log_{\sigma}(\bepsilon/2),
\end{equation*}
which together with $0<\sigma<1$ implies that $\sigma^\bK \leq \bepsilon/2$. It follows from this,  $0<\eta_0\leq1$, and $\eta_\bK = \eta_0 \sigma^\bK$  that $\eta_\bK\leq\bepsilon/2$. This relation and \eqref{relative-change} lead to a contradiction to \eqref{non-termination}. Hence, Algorithm \ref{finitesum:alg1} terminates in at most $\bK+1$ outer iterations. 

Suppose now that the termination condition $\rho_k^{-1} \left\| (x_{k+1}, \lambda_{k+1}) - (x_k, \lambda_k) \right\| +\eta_k \leq \bepsilon$  holds at some iterate $k$ for Algorithm \ref{finitesum:alg1}. Then, by an argument similar to that in the proof of \cite[Theorem 4(i)]{lu2023iteration}, it follows that the output $(x_{k+1},\lambda_{k+1})$ is an $\bepsilon$-KKT solution of problem \eqref{primal} and its dual.

In addition, note that Algorithm~\ref{alg4} is called at the $k$th outer iteration of Algorithm~\ref{finitesum:alg1} to solve the proximal AL subproblem~\eqref{PAL-subprob}. Moreover, each call of Algorithm~\ref{alg4} performs the same number of evaluations of $\nabla \barf$ and $\nabla c$. Let $n_k$ denote the number of evaluations of either $\nabla \barf$ or $\nabla c$ performed by Algorithm~\ref{alg4} during the $k$th outer iteration of Algorithm~\ref{finitesum:alg1}. Note that when applied to~\eqref{PAL-subprob}, Algorithm~\ref{alg4} performs at most two evaluations of $\nabla \barf$ and $\nabla c$ per iteration. By this fact, Lemma \ref{lem:Lk}, and Theorem \ref{thm:cplx-fom-str-cvx} with $\tf(\cdot) \gets \overline{\cL}(\cdot, \lambda_k;\rho_k)-\psi(\cdot)$, $L_{\nabla \tf} \gets L_k$, $\mu_\tf \gets \mu+\rho^{-1}_k$, and $\varepsilon \gets \eta_k$, one has
\beq \label{nk}
n_k \leq 2\left\lceil \sqrt{\frac{L_k}{\mu+\rho^{-1}_k}}\ \right\rceil \max\left\{ 1, \left\lceil 2\log \frac{2L_kD_X}{\eta_k}\right\rceil\right\},
\eeq
where $L_k$ is given in \eqref{L_k}.

Further, using $\rho_0\geq1$, $\zeta>1$, $0<\zeta\sigma <1$, $\rho_i= \rho_0 \zeta^i$, $\eta_ i= \eta_0 \sigma^i$, \eqref{alg1-para}, and Lemma \ref{lem:Lk}, we obtain that
\begin{align}
L_k  & \leq  C \rho_k+B+L_{\nabla c}\sum^{k-1}_{i=0}\rho_i\eta_i+1 = C \rho_0 \zeta^k+B+L_{\nabla c}\rho_0\eta_0\sum^{k-1}_{i=0} (\zeta\sigma)^i+1 \nn \\
& \leq C \rho_0 \zeta^k+B+\frac{L_{\nabla c}\rho_0\eta_0}{1-\zeta\sigma}+1  \overset{\eqref{alg1-para}}{\leq} \bar{C}\zeta^k. \label{Lk_bound}
\end{align}
Notice that $L_{\nabla \barf} \geq \mu$. By this fact and \eqref{alg1-para}, one can observe that $\bar{C} \geq \mu$, which together with $0<\eta_0, \sigma \leq 1$ and $\bar{C}, \zeta, D_X \geq 1$ implies that 
\beq \label{aux-ineq}
\sqrt{{\bar{C}}/{\mu}}\,\zeta^{{k}/{2}}\geq1, \qquad 2 \log \big({2\bar{C}D_X\zeta^k }/({\eta_0\sigma^k})\big)\geq2\log2\geq1. 
\eeq 
Using these relations, \eqref{nk}, \eqref{Lk_bound}, and $\rho_k \geq 1$,  we obtain that 
\begin{align}
    n_k&\overset{\eqref{nk}}{\leq} 2 \left\lceil \sqrt{\frac{L_k}{\mu}}\, \right\rceil\max \left\{ 1,\; \left\lceil 2 \log \frac{2  L_k D_X}{\eta_k} \right\rceil \right\} 
   \overset{\eqref{Lk_bound}}{ \leq}2 \left\lceil \sqrt{\frac{\bar{C}}{\mu}}\zeta^{\frac{k}{2}} \right\rceil\max \left\{ 1,\; \left\lceil 2 \log \frac{2 \bar{C} D_X\zeta^k }{\eta_0\sigma^k} \right\rceil \right\} \nn \\
    &\leq 2 \left\lceil \sqrt{\frac{\bar{C}}{\mu}}\zeta^{\frac{k}{2}} \right\rceil \left\lceil 2 \log \frac{2 \bar{C} D_X\zeta^k}{\eta_0\sigma^k}  \right\rceil  
    \overset{\eqref{aux-ineq}}{\leq} 16\sqrt{\frac{\bar{C}}{\mu}}\zeta^{\frac{k}{2}}   \log \frac{2\bar{C}D_X\zeta^k}{\eta_0\sigma^k}\nn \\ 
& \leq 16\sqrt{\frac{\bar{C}}{\mu}}   \log \left( \frac{2 \bar{C} D_X\zeta }{\eta_0\sigma}\right) 
\max\big\{1, k \zeta^{k/2}\big\}, \label{nk-bound}
\end{align}
where the last inequality follows from $2\bar{C} D_X\zeta /\eta_0 \geq 1$ due to $\bar{C},  D_X, \zeta \geq 1$,  and $0<\eta_0\leq1$. Observe that
\begin{align}
1+\sum_{k=1}^\bK k\zeta^{k/2} &\leq  \bK \sum_{k=0}^\bK \zeta^{k/2} \leq  \bK \zeta^{(\bK+1)/2}/(\sqrt{\zeta}-1),  \label{sum-zetak1}\\
\zeta^{(\bK+1)/2} & \overset{\eqref{bK}}{\leq}  \zeta \max\left\{1, \sqrt{\frac{2(D+\rho_0\eta_0)}{(1-\zeta\sigma)\epsilon}}\right\}. \label{sum-zetak2}
\end{align}
By these relations, \eqref{bK},  \eqref{bN}, and \eqref{nk-bound}, one has
\begin{align*}
\sum_{k=0}^{\bK} n_k &=n_0+\sum_{k=1}^{\bK} n_k \overset{\eqref{nk-bound}}{\leq} 16 \sqrt{\frac{\bar{C}}{\mu}}\log\Big(\frac{2 \bar{C} D_X\zeta}{\eta_0\sigma}\Big) \Big(1+\sum_{k=1}^{\bK} k\zeta^{\frac{k}{2}}\Big) \\
& \overset{\eqref{sum-zetak1}}{\leq} \frac{16\sqrt{\bar{C}}}{\sqrt{\mu}(\sqrt{\zeta}-1)}\log\Big(\frac{2\bar{C} D_X\zeta}{\eta_0\sigma}\Big)\bK \zeta^{(\bK+1)/2} \overset{\eqref{bK}\eqref{bN}\eqref{sum-zetak2}}{\leq} \bN,
\end{align*}
and hence the conclusion of Theorem \ref{thm:sub-complexity} holds.
\end{proof}

\subsection{Proof of the main results in Section \ref{sec:finite-sum}} \label{proof:thm2}

{ In this subsection, we present a proof of Theorem~\ref{thm:finite-sum}, which establishes the convergence properties of Algorithm~\ref{alg2}—an inexact accelerated SVRG method for problem~\eqref{prob1}.

Before proceeding, we highlight the main technical challenges in analyzing Algorithm~\ref{alg2} compared to the classical (exact) accelerated SVRG method~\cite{lan2019unified}. To this end, we introduce the following notation:
\begin{align}
 &\ell_{f}(u, v):=f(v)+\langle \nabla f(v), u-v\rangle, \qquad \delta_t:=g_t-\nabla f(y_t), \label{notation} \\
 & \Omega:=\{x\in X: c(x) \in -\cK\}, \qquad
 \Omega_\epsilon :=\{x\in X:\dist(c(x),-\cK)\leq\epsilon\}.  \label{omega-eps}
\end{align}

In the classical accelerated SVRG method, $z_t$ is assumed to be the exact optimal solution of the subproblem~\eqref{subproblem} in Algorithm~\ref{alg2}. Under this assumption, the following key inequality holds:
\beq \label{exact-svrg}
 \gamma_k\big(\ell_{f}(z_t, y_t)-\ell_{f}(x,y_t)+\psi(z_t)-\psi(x)\big)\leq \frac{1}{2}\|x-z_{t-1}\|^2-\frac{1}{2}\|x-z_{t}\|^2-\frac{1}{2}\|z_{t-1}-z_t\|^2+\gamma_k\langle \delta_t, x-z_t\rangle 
\eeq
for all $1\leq t\leq T_k$, $k \geq 1$, and $x\in\Omega$ (see, e.g.,~\cite[Lemma~5]{lan2019unified}). The analysis of the classical accelerated SVRG method relies crucially on this inequality.

However, in Algorithm~\ref{alg2}, $z_t$ is only an approximate solution of the subproblem~\eqref{subproblem}, and thus~\eqref{exact-svrg} does not generally hold. Instead, by leveraging the fact that $(z_t,\lambda_t)$ is an $\epsilon$-KKT pair for~\eqref{subproblem}, together with several technical arguments, we establish the following inequality:
\[
 \gamma_k\big(\ell_{f}(z_t, y_t)-\ell_{f}(x,y_t)+\psi(z_t)-\psi(x)\big)\leq \frac{1}{2}\|x-z_{t-1}\|^2-\frac{1}{2}\|x-z_{t}\|^2-\frac{1}{2}\|z_{t-1}-z_t\|^2+\gamma_k\langle \delta_t, x-z_t\rangle+\cO(\gamma_k\epsilon) 
\]
for all $1\leq t\leq T_k$, $k \geq 1$, and $x\in\Omega_\epsilon$, where $\{\gamma_k\}$ is the stepsize sequence of Algorithm~\ref{alg2} (see Lemma~\ref{l-eoptimality} for details). Using this inequality, together with an appropriate choice of $\{\gamma_k\}$, we show that the outer loop of Algorithm~\ref{alg2} achieves the same convergence rate as the classical accelerated SVRG method, despite solving its subproblems inexactly.
}

We now present the proof of Theorem~\ref{thm:finite-sum}. To this end, we begin by establishing several technical lemmas. The following lemma provides a uniform upper bound on the optimal Lagrange multipliers of problem~\eqref{subproblem}.

\begin{lemma}\label{l-subbound}
Suppose that Assumption \ref{a1} holds.  Let $\Lambda$ be defined in \eqref{Lambda}, and $\lambda_{k,t}^*$ be an optimal Lagrange multiplier of problem \eqref{subproblem} that is generated in Algorithm \ref{alg2} with $\{\gamma_k\}$ and $\{q_i\}$ given in \eqref{def-para}.\footnote{The existence of $\lambda_{k,t}^*$ is guaranteed by the strong convexity of problem \eqref{subproblem} and Assumption \ref{a1}(iv).} Then we have  
\beq \label{lambda*}
  \|\lambda_{k,t}^*\|\leq\Lambda \qquad\forall  1 \leq t \leq T_k, k \geq 1.
\eeq
\end{lemma}

\begin{proof}
Fix any $k \geq 1$ and $1 \leq t \leq T_k$. Observe from Algorithm \ref{alg2} that $y_t, \tilde x_{k}\in X$, and 
\[
g_t=(\nabla f_{i_t}(y_t)-\nabla f_{i_t}(\tilde x_{k}))/(q_{i_t}s)+\tilde g_{k}.
\]
Using these observations,  \eqref{DX},  \eqref{DF}, \eqref{Lf}, \eqref{f-bnd}, \eqref{def-para}, and Assumption \ref{a1}, we obtain that
\begin{align}
\|g_t\|&=\left\|(\nabla f_{i_t}(y_t)-\nabla f_{i_t}(\tilde x_{k}))/(q_{i_t}s)+\tilde g_{k}\right\|\leq L_{\nabla f_{i_t}}\|y_t-\tilde x_k\|/(q_{i_t}s)+\|\tilde g_k\|\nn\\
&\overset{\eqref{Lf}\eqref{def-para}}{=}L_{\nabla f}\|y_t-\tilde x_k\|+\|\nabla f(\tilde x_k)\|
\overset{\eqref{DX}\eqref{f-bnd}}{\leq} L_{\nabla f}D_X+L_f.\label{gt}
\end{align}
Notice from \eqref{def-para} that $\gamma_k^{-1}=3L_{\nabla f}\alpha_k\leq3L_{\nabla f}/2$. By this relation,  \eqref{DX},  \eqref{DF},  \eqref{gt}, and the definition of $F_{k,t}$ in \eqref{subproblem}, one has that for all $x, y\in X$, 
\begin{align}
F_{k,t}(x)-F_{k,t}(y)&
\overset{\eqref{subproblem}}{=}\frac{1}{2\gamma_k}\big(\|x-z_{t-1}+\gamma_kg_t\|^2-\|y-z_{t-1}+\gamma_kg_t\|^2\big)+\psi(x)-\psi(y)\nn\\
&=\frac{1}{2\gamma_k}\langle x-y, x-z_{t-1}+y-z_{t-1}+2\gamma_kg_t\rangle +\psi(x)-\psi(y)\nn\\
&\leq\frac{1}{2\gamma_k}\|x-y\| \|x-z_{t-1}+y-z_{t-1}+2\gamma_kg_t\| +\psi(x)-\psi(y)\nn\\
&\leq\frac{1}{2\gamma_k}\|x-y\| \big(\|x-z_{t-1}\|+\|y-z_{t-1}\|+2\gamma_k\|g_t\|\big)+\psi(x)-\psi(y)\nn\\
&\leq \frac{1}{2\gamma_k}D_X\big(2D_X+2\gamma_k(L_{\nabla f}D_X+L_f)\big)+D_\psi\nn\\
&=D_X^2 \gamma_k^{-1}+(L_{\nabla f}D_X+L_f)D_X+D_\psi
\leq 5L_{\nabla f}D_X^2/2+L_f D_X+D_\psi, \label{Fktbound}
\end{align}
where the third inequality follows from \eqref{DX}, \eqref{DF}, and \eqref{gt}, and the last inequality is due to 
$\gamma_k^{-1}\leq3L_{\nabla f}/2$.

Let $x_{k,t}^*$ be the optimal solution of problem \eqref{subproblem}. Since $\lambda_{k,t}^*$ is an optimal Lagrange multiplier of \eqref{subproblem}, one has that $\lambda_{k,t}^*\in\cK^*$,  $\langle\lambda_{k,t}^*,c(x_{k,t}^*)\rangle=0$, and 
\[
x_{k,t}^* = \argmin_x \big\{F_{k,t}(x)+\langle\lambda_{k,t}^*,c(x)\rangle \big\}. 
\]
Notice from Assumption \ref{a1}(iv) that $\xint \in X$ and $c(\xint)+D_c\cB \subset -\cK$. It then follows that
\[
F_{k,t}(\xint)+\langle\lambda_{k,t}^*,c(\xint)\rangle \geq F_{k,t}(x_{k,t}^*)+\langle\lambda_{k,t}^*,c(x_{k,t}^*)\rangle= F_{k,t}(x_{k,t}^*),
\]
which together with  \eqref{Fktbound} implies that
\begin{equation}\label{l1-1}
\langle\lambda_{k,t}^*,-c(\xint)\rangle\leq F_{k,t}(\xint)-F_{k,t}(x_{k,t}^*)
\leq 5L_{\nabla f}D_X^2/2+L_f D_X+D_\psi.  
\end{equation}
If $\lambda_{k,t}^*=0$, \eqref{lambda*} clearly holds. Now suppose $\lambda_{k,t}^*\neq 0$. 
By $c(\xint)+D_c\cB \subset -\cK$, one has $-c(\xint)-D_c\lambda_{k,t}^*/\|\lambda_{k,t}^*\|\in\cK$. Using this and $\lambda_{k,t}^*\in \cK^*$, we obtain that
\begin{equation*}
\langle\lambda_{k,t}^*,-c(\xint)-D_c\lambda_{k,t}^*/\|\lambda_{k,t}^*\|\rangle\geq0,
\end{equation*}
which implies that $\|\lambda_{k,t}^*\|\leq D_c^{-1} \langle\lambda_{k,t}^*,-c(\xint)\rangle$. Combining this relation with \eqref{Lambda} and \eqref{l1-1}, we conclude that \eqref{lambda*} holds.
\end{proof}

The following lemma establishes the convexity of the set $\Omega_\epsilon$ defined in~\eqref{omega-eps}.

\begin{lemma}\label{l-feasibility}
Let $\Omega_\epsilon$ be defined in \eqref{omega-eps}. Then $\Omega_\epsilon$ is convex for any $\epsilon>0$.
\end{lemma}

\begin{proof}
Fix any $x,y\in \Omega_\epsilon$ and $\alpha\in[0,1]$. In view of \eqref{K-convex} and the fact that $\cK$ is a convex cone, one has 
\[
\alpha \cK + (1-\alpha)\cK=\cK, \qquad -\cK+c(\alpha x+(1-\alpha)y)-\alpha c(x)-(1-\alpha)c(y) \subseteq -\cK.
\]
Using these relations, we obtain that
\begin{align*}
&\alpha\,\dist\big(c(x),-\cK\big)+(1-\alpha)\dist\big(c(y),-\cK\big) 
 \geq \dist\big(\alpha c(x)+(1-\alpha)c(y),-\cK\big) \\ 
 &= \dist\big(c(\alpha x+(1-\alpha)y),-\cK+c(\alpha x+(1-\alpha)y)-\alpha c(x)-(1-\alpha)c(y)\big) \\
& \geq \dist\big(c(\alpha x+(1-\alpha)y),-\cK\big),
\end{align*}
which implies that $\dist(c(\cdot),-\cK)$ is a convex function.  The convexity of $\Omega_\epsilon$ then follows from this fact and the convexity of $X$.
\end{proof}

The next lemma establishes some properties of the sequences $\{\tilde x_k\}$, $\{\tilde z_k\}$, $\{x_t\}$,  $\{y_{t}\}$, $\{z_t\}$, and $\{\lambda_t\}$ generated by Algorithm \ref{alg2}, which will be used in the subsequent analysis. 

\begin{lemma} \label{iter-property}
Suppose that Assumption \ref{a1} holds. Let $D_X$, $\Lambda$,  $\Omega_\epsilon$, $\rho_0$, $\eta_0$, $\zeta$, $\sigma$, $R$ be given in \eqref{DX}, \eqref{Lambda}, \eqref{omega-eps}, and Algorithms \ref{finitesum:alg1} and \ref{alg2}, respectively, and $\{\tilde x_k\}$, $\{\tilde z_k\}$, $\{x_t\}$, $\{y_t\}$, $\{z_t\}$, and $\{\lambda_t\}$ be generated in Algorithm \ref{alg2}. Then the following statements hold.
\bi
\item[(i)] $\tilde x_k, \tilde z_k \in \Omega_\epsilon$ for all $k \geq 1$.
\item[(ii)] $x_t, z_t \in \Omega_\epsilon$ for all $0 \leq t \leq T_k$ and $k \geq 1$, and $y_t \in \Omega_\epsilon$ for all $1 \leq t \leq T_k$ and $k \geq 1$.
\item[(iii)] $\|(x_t,\lambda_t)-(x^*_{k,t},\lambda^*_{k,t})\| \leq  D_X + R+\Lambda+\rho_0\eta_0/(1-\zeta\sigma)$ for all $1 \leq t \leq T_k$ and $k \geq 1$, where $(x^*_{k,t}, \lambda^*_{k,t})$ is an arbitrary pair of optimal solution and Lagrange multiplier of problem \eqref{subproblem}.
\ei
\end{lemma}

\begin{proof}
For all $1 \leq t \leq T_k$ and $k \geq 1$, it follows from Theorem \ref{thm:sub-complexity} that $(z_t,\lambda_t)$ is a pair of $\epsilon$-KKT solution of problem \eqref{subproblem}. By this fact and the definition of $\epsilon$-KKT solution in \eqref{ekkt}, one can see that $z_t\in X$, $\lambda_t \in \cK^*$, and $\dist(c(z_t),\cN_{\cK^*}(\lambda_t)) \leq \epsilon$. Hence, there exists $u\in \cN_{\cK^*}(\lambda_t)$  such that $\|c(z_t)-u\| \leq \epsilon$. By $\lambda_t \in \cK^*$, $u\in \cN_{\cK^*}(\lambda_t)$, and the fact that 
$\cK^*$ is a cone, one can observe that $u^T \lambda_t=0$ and $u^T(v- \lambda_t) \leq 0$ for all $v\in\cK^*$. It follows that $u^Tv \leq 0$ for all $v\in\cK^*$, which leads to $u\in -\cK$. Using this and $\|c(z_t)-u\| \leq \epsilon$,  we see that $\dist(c(z_t),-\cK) \leq \epsilon$, which together with $z_t\in X$ implies that $z_t\in \Omega_\epsilon$ for all $1 \leq t \leq T_k$ and $k \geq 1$. In addition, notice from Algorithm \ref{alg2} that $x^0\in \Omega_\epsilon$ and $\tx_1=\tz_1=x^0$.  Hence, $\tilde x_1, \tilde z_1 \in \Omega_\epsilon$. Also, for each $k\geq 1$, $x_0=\tilde x_{k}$, $z_0=\tilde z_{k}$, $\tilde z_{k+1}=z_{T_k}$, $\tilde x_{k+1}$ is a convex combination of $x_1, \ldots, x_{T_k}$,  $y_t$ is a convex combination of $x_{t-1}, z_{t-1}, \tilde x_{k}$, and $x_t$ is a convex combination of $x_{t-1}$,  $z_{t}$, and $\tilde x_{k}$  for all $1 \leq t \leq T_k$. Based on these facts, the convexity of $\Omega_\epsilon$, and an inductive argument, one can readily verify that statements (i) and (ii) hold.

We next show that statement (iii) holds. Indeed, since $(x_t,\lambda_t)$ is generated by Algorithm \ref{finitesum:alg1} applied to problem \eqref{subproblem}, it follows from \eqref{iter-bound} that 
\beq \label{opt-gap}
\|(x_t,\lambda_t)-(x^*_{k,t},\lambda^*_{k,t})\| \leq \|(x_0,\lambda_0)-(x^*_{k,t},\lambda^*_{k,t})\| +\sum^\infty_{i=0}\rho_i\eta_i \qquad \forall 0 \leq t \leq T_k, k \geq 1.
\eeq
Notice from Algorithm \ref{alg2} that $x_0, x^*_{k,t}\in X$, $\|\lambda_0\| \leq R$,  $\rho_i= \rho_0 \zeta^i$, and $\eta_i = \eta_0 \sigma^i$ for all $i \geq 0$. Also, recall from Lemma 
\ref{l-subbound} that $\|\lambda^*_{k,t}\| \leq \Lambda$. In view of these facts, Assumption \ref{a1}, and $0<\zeta \sigma<1$, one has
\begin{align*}
 \|(x_0,\lambda_0)-(x^*_{k,t},\lambda^*_{k,t})\| +\sum^\infty_{i=0}\rho_i\eta_i  
&\leq \|x_0-x^*_{k,t}\|+\|\lambda_0\|+\|\lambda^*_{k,t}\| + \rho_0\eta_0\sum^\infty_{i=0}(\zeta\sigma)^i \\
& \leq D_X + R + \Lambda+\rho_0\eta_0/(1-\zeta\sigma),
\end{align*}
which together with \eqref{opt-gap} implies that statement (iii) holds as desired.
\end{proof} 

The following lemma establishes a key property of the inexact solution $z_t$ to problem \eqref{subproblem} generated in Algorithm \ref{alg2}.

\begin{lemma}\label{l-eoptimality}
Suppose that Assumption \ref{a1} holds.  Let $f$, $C_5$, $\ell_f$, $\delta_t$, $\Omega_\epsilon$, and $\epsilon$ be given in \eqref{Lf},  \eqref{c45}, \eqref{notation}, \eqref{omega-eps},  and Algorithm \ref{alg2}, respectively, and let $\{g_t\}$, $\{y_t\}$, $\{z_t\}$ be generated in the $k$th outer iteration of Algorithm \ref{alg2} with $\{\gamma_k\}$ and $\{q_i\}$ given in \eqref{def-para}. Then,  for all $1\leq t\leq T_k$, $k \geq 1$, and $x\in\Omega_\epsilon$, we have
\begin{equation}\label{eoptimality}
    \gamma_k\big(\ell_{f}(z_t, y_t)-\ell_{f}(x,y_t)+\psi(z_t)-\psi(x)\big)\leq \frac{1}{2}\|x-z_{t-1}\|^2-\frac{1}{2}\|x-z_{t}\|^2-\frac{1}{2}\|z_{t-1}-z_t\|^2+\gamma_k\langle \delta_t, x-z_t\rangle+\gamma_kC_5 \epsilon.
\end{equation}
\end{lemma}

\begin{proof}
Fix any $k \geq 1$. Notice from Algorithm \ref{alg2} that $(z_t,\lambda_t)$ is a pair of $\epsilon$-KKT solution of problem \eqref{subproblem}. By this fact and the definition of $\epsilon$-KKT solution in \eqref{ekkt}, one has that $z_t\in X$, $\lambda_t \in \cK^*$,  and
\begin{equation}\label{eopt-1}
\dist\big(0,(z_t-z_{t-1}+\gamma_kg_t)/\gamma_k+\nabla c(z_t)\lambda_t+\partial\psi(z_t)\big)\leq \epsilon
, \quad \dist\big(c(z_t),\cN_{\cK^*}(\lambda_t)\big)\leq \epsilon.
\end{equation}
Hence,  there exists $u\in \cN_{\cK^*}(\lambda_t)$ such that $\|c(z_t)-u\|\leq \epsilon$. By $\lambda_t \in \cK^*$, $u\in \cN_{\cK^*}(\lambda_t)$, and the fact that 
$\cK^*$ is a cone, one can observe that $\langle \lambda_t, u\rangle=0$. In addition, let $\lambda_{k,t}^*$ be an optimal Lagrange multiplier of problem \eqref{subproblem}. Invoking Lemmas \ref{l-subbound} and \ref{iter-property}(iii), we have
\[
\|\lambda_{k,t}^*\| \leq \Lambda, \qquad \|\lambda_t-\lambda_{k,t}^*\| \leq D_X + R+\Lambda+\rho_0\eta_0/(1-\zeta\sigma),
\]
which imply that 
\beq \label{lambda-bound}
\|\lambda_t\| \leq \|\lambda_{k,t}^*\|+\|\lambda_t-\lambda_{k,t}^*\| \leq D_X + R+2\Lambda+\rho_0\eta_0/(1-\zeta\sigma).
\eeq
Using this relation,  $\langle \lambda_t, u\rangle=0$, and $\|c(z_t)-u\|\leq \epsilon$, we obtain that
\begin{align}\label{eopt-2}
 \langle\lambda_t,c(z_t)\rangle=\langle\lambda_t,c(z_t)-u\rangle \geq-\epsilon\|\lambda_t\| \geq-\epsilon\big(D_X + R+2\Lambda+\rho_0\eta_0/(1-\zeta\sigma)\big).
\end{align}
Similarly,  for any $x\in\Omega_\epsilon$, it follows from the second relation in \eqref{omega-eps} that there exists $\tu \in-\cK$ such that $\|c(x)-\tu\| \leq \epsilon$. This  together with $\lambda_t\in\cK^*$ and \eqref{lambda-bound} implies 
that $\langle\lambda_t,\tu\rangle \leq 0$, and hence
\beq  \label{eopt-3}
 \langle\lambda_t,c(x)\rangle  \leq\langle\lambda_t,c(x)-\tu\rangle 
\leq \epsilon\|\lambda_t\| \overset{\eqref{lambda-bound}}{\leq}\epsilon\big(D_X + R+2\Lambda+\rho_0\eta_0/(1-\zeta\sigma)\big) \qquad \forall x\in\Omega_\epsilon.
\eeq
Also, by the first relation in \eqref{eopt-1}, there exists a vector $\psi'(z_t)\in\partial\psi(z_t)$ such that 
\beq \label{eps-kkt}
\|(z_t-z_{t-1}+\gamma_kg_t)/\gamma_k+\nabla c(z_t)\lambda_t+\psi'(z_t)\|\leq\epsilon.
\eeq
In addition, by the $\cK$-convexity of $c(\cdot)$ and $\lambda_t\in\cK^*$, one has
\begin{equation*}
\langle\lambda_t,\alpha c(x)+(1-\alpha)c(y)-c(\alpha x+(1-\alpha)y)\rangle\geq0\qquad\forall x,y\in X,\alpha\in[0,1],
\end{equation*}
which implies that the function $\langle\lambda_t,c(\cdot)\rangle$ is convex. It then follows that the function 
$\|x-z_{t-1}+\gamma_kg_t\|^2/(2\gamma_k){ +\langle \lambda_t,c(x)\rangle}+\psi(x)$ 
is strongly convex with modulus $\gamma_k^{-1}$.  Using this fact,  \eqref{eps-kkt}, and Assumption \ref{a1}, we obtain that for all $x\in \Omega_\epsilon$,
\begin{align*}
&\frac{\|x-z_{t-1}+\gamma_kg_t\|^2}{2\gamma_k}+\langle\lambda_t,c(x)\rangle +\psi(x)\geq \frac{\|z_t-z_{t-1}+\gamma_kg_t\|^2}{2\gamma_k}+\langle\lambda_t,c(z_t)\rangle+\psi(z_t)\nn\\
&\quad\ +\langle(z_t-z_{t-1}+\gamma_kg_t)/\gamma_k+\nabla c(z_t)\lambda_t+\psi'(z_t),x-z_t\rangle+\frac{1}{2\gamma_k}\|x-z_t\|^2\nn\\
&\overset{\eqref{eps-kkt}}{\geq}\frac{\|z_t-z_{t-1}+\gamma_kg_t\|^2}{2\gamma_k}+\langle\lambda_t,c(z_t)\rangle+\psi(z_t)-\|x-z_t\|\epsilon+\frac{1}{2\gamma_k}\|x-z_t\|^2\nn\\
&\geq\frac{\|z_t-z_{t-1}+\gamma_kg_t\|^2}{2\gamma_k}+\langle\lambda_t,c(z_t)\rangle+\psi(z_t)- D_X \epsilon+\frac{1}{2\gamma_k}\|x-z_t\|^2,
\end{align*}
where the last inequality follows from the fact that $\|x-z_t\| \leq D_X$ for all $x\in \Omega_\epsilon$ due to Assumption \ref{a1}.  Combining the above inequality with \eqref{eopt-2} and \eqref{eopt-3} leads to
\begin{align}
\frac{\|x-z_{t-1}+\gamma_kg_t\|^2}{2\gamma_k}+\psi(x) & \geq\frac{\|z_t-z_{t-1}+\gamma_kg_t\|^2}{2\gamma_k}+\psi(z_t)+\frac{1}{2\gamma_k}\|x-z_t\|^2 \nn \\
& \quad-\big(3D_X + 2R+4\Lambda+2\rho_0\eta_0/(1-\zeta\sigma)\big)\epsilon \qquad \forall x\in \Omega_\epsilon.\label{eopt-4}
\end{align}
Moreover, by  \eqref{notation}, one has
\begin{align*}
& \|z_t-z_{t-1}+\gamma_kg_t\|^2-\|x-z_{t-1}+\gamma_kg_t\|^2=\|z_t-z_{t-1}\|^2-\|x-z_{t-1}\|^2+2\gamma_k\langle g_t,z_t-x\rangle\\
&\overset{\eqref{notation}}{=}\|z_t-z_{t-1}\|^2-\|x-z_{t-1}\|^2+2\gamma_k\langle \delta_t,z_t-x\rangle+2\gamma_k\langle \nabla f(y_t),z_t-x\rangle\\
&\overset{\eqref{notation}}{=}\|z_t-z_{t-1}\|^2-\|x-z_{t-1}\|^2+2\gamma_k\langle \delta_t,z_t-x\rangle+2\gamma_k\big(\ell_{f}(z_t, y_t)-\ell_{f}(x,y_t )\big).
\end{align*}
{ Substituting this relation into~\eqref{eopt-4} and using~\eqref{c45}, we obtain \eqref{eoptimality}, which completes the proof of the lemma.}
\end{proof}

The next lemma, established in \cite[Lemma 3]{lan2019unified}, shows that $g_t$ is an unbiased stochastic estimator of $\nabla f(y_t)$ and provides an upper bound on its variance.

\begin{lemma}[\textbf{\cite[Lemma 3]{lan2019unified}}]
Suppose that Assumption \ref{a1} holds. Let $f$, $L_{\nabla f}$  and $\delta_t$ be defined in \eqref{Lf} and \eqref{notation}, respectively, and $\tilde x_{k}$ and $\{y_t\}$ be generated in the $k$th outer iteration of Algorithm \ref{alg2} with $\{q_i\}$ given in \eqref{def-para}.
Then,  for all $1\leq t\leq T_k$ and $k \geq 1$, we have
\beq
    \bbE[\delta_t|\Xi_{k,t-1}]=0, \qquad
    \bbE[\|\delta_t\|^2|\Xi_{k,t-1}]\leq 2L_{\nabla f}\big(f(\tilde x_{k}) - f(y_t)-\langle\nabla f(y_t),\tilde x_{k} - y_t \rangle\big),\label{mean-var}
\eeq
where $\Xi_{k,t-1}$ denotes the collection of all samples generated by Algorithm~\ref{alg2} prior to the $t$-th inner iteration of the $k$th outer iteration.
\end{lemma}

The following lemma establishes a recursion for $\{\bbE[F(x_t)-F(x)]\}$ for all $1\leq t\leq T_k$, $k \geq 1$, and $x\in\Omega_\epsilon$.

\begin{lemma}\label{l-iter2}
Suppose that Assumption \ref{a1} holds. Let $L_{\nabla f}$, $C_5$, $\Omega_\epsilon$, $\epsilon$ be given in  \eqref{Lf}, \eqref{c45}, \eqref{omega-eps}, and Algorithm \ref{alg2}, respectively, and let $\{x_t\}$, $\{y_t\}$, $\{z_t\}$,  $\{g_t\}$, and $\tilde x_{k}$ be generated in the $k$th outer iteration of Algorithm \ref{alg2} with $\{q_i\}$ and $\gamma_k$ given in \eqref{def-para}, and $\alpha_k$, $\gamma_k$, and $p_k$ satisfying
\begin{align}
    &1-L_{\nabla f}\alpha_k\gamma_k>0, \label{l7-1}\\
    &p_k-\frac{L_{\nabla f}\alpha_k \gamma_k}{1-L_{\nabla f} \alpha_k\gamma_k}\geq 0 \label{l7-2}
\end{align}
for all $k\geq1$. Then,  for all $1\leq t\leq T_k$, $k \geq 1$, and $x\in\Omega_\epsilon$, we have
\begin{align}
    \frac{\gamma_k}{\alpha_k}\bbE[F(x_t)-F(x)]+\frac{1}{2}\bbE[\|x-z_t\|^2] &\leq \frac{\gamma_k}{\alpha_k}(1-\alpha_k-p_k)\bbE[F(x_{t-1})-F(x)]+\frac{\gamma_kp_k}{\alpha_k} \bbE[F(\tilde x_{k})-F(x)]\nn\\
    &\quad+\frac{1}{2}\bbE[\|x-z_{t-1}\|^2]+\gamma_kC_5 \epsilon. \label{l7-exp}
\end{align}
\end{lemma}

\begin{proof}
Fix any $k \geq 1$.  Let $\{x_t\}$, $\{y_t\}$, and $\{z_t\}$ be generated in the $k$th outer iteration of Algorithm \ref{alg2} for all $1\leq t\leq T_k$. Then it follows from the expressions of $x_t$ and $y_t$ in Algorithm \ref{alg2} that $ x_t-y_t=\alpha_k(z_t-z_{t-1})$.
By this relation, \eqref{eoptimality}, the $L_{\nabla f}$-smoothness and the convexity of $f$, and the expressions of $x_t$ and $\ell_{f}(\cdot,\cdot)$, one has
\begin{align}
    &f(x_t)\leq \ell_{f}(x_t,y_t)+\frac{L_{\nabla f}}{2}\|x_t-y_t\|^2\nn\\
    &=(1-\alpha_k-p_k)\ell_{f}(x_{t-1},y_t)+\alpha_k \ell_{f}(z_t,y_t)+p_k\ell_{f}(\tilde x_{k},y_t)+\frac{L_{\nabla f}\alpha_k^2}{2}\|z_t-z_{t-1}\|^2\nn\\
    &\overset{\eqref{eoptimality}}{\leq}(1-\alpha_k-p_k)\ell_{f}(x_{t-1},y_t) + \alpha_k\Big(\ell_{f}(x,y_t)+\psi(x)-\psi(z_t)+\frac{1}{2\gamma_k}\big(\|x-z_{t-1}\|^2-\|x-z_t\|^2\nn \\ 
&\qquad-\|z_{t-1}-z_t\|^2\big)+\langle \delta_t, x-z_t\rangle+C_5\epsilon\Big)
+p_k\ell_{f}(\tilde x_{k},y_t)+\frac{L_{\nabla f}\alpha_k^2}{2}\|z_t-z_{t-1}\|^2\nn\\
    &= (1-\alpha_k-p_k)\ell_{f}(x_{t-1},y_t)
    +  \alpha_k\Big(\ell_{f}(x,y_t)+\psi(x)-\psi(z_t)+\frac{1}{2\gamma_k}\|x-z_{t-1}\|^2-\frac{1}{2\gamma_k}\|x-z_t\|^2\Big)\nn\\
    &\quad +p_k\ell_{f}(\tilde x_{k},y_t)-\frac{\alpha_k}{2\gamma_k}(1-L_{\nabla f}\alpha_k\gamma_k)\|z_{t-1}-z_t\|^2+\alpha_k\langle \delta_t, x-z_t\rangle+\alpha_kC_5\epsilon \nn\\
    &\leq(1-\alpha_k-p_k)f(x_{t-1})
    +\alpha_k\Big(F(x)-\psi(z_t)+\frac{1}{2\gamma_k}\|x-z_{t-1}\|^2-\frac{1}{2\gamma_k}\|x-z_t\|^2\Big)\nn\\
    &\quad + p_k\ell_{f}(\tilde x_{k},y_t)-\frac{\alpha_k}{2\gamma_k}(1-L_{\nabla f}\alpha_k\gamma_k)\|z_{t-1}-z_t\|^2+\alpha_k\langle \delta_t, x-z_t\rangle+\alpha_kC_5\epsilon, \nn
\end{align}
where the first relation follows from the $L_{\nabla f}$-smoothness of $f$, the second relation is due to the definition of $\ell_{f}(\cdot,\cdot)$, and the last relation follows from the convexity of $f$ and the expression of $F$.  Also, by the expression of $x_t$ and the convexity of $\psi$, one has
\[
    \psi(x_t)\leq (1-\alpha_k-p_k)\psi(x_{t-1})+\alpha_k\psi(z_t)+p_k\psi(\tilde x_{k}).
\]
Summing up these two inequalities and using the expression of $F$, we have
\begin{align}
    F(x_t)& \leq(1-\alpha_k-p_k)F(x_{t-1})+\alpha_k\Big(F(x)+\frac{1}{2\gamma_k}\|x-z_{t-1}\|^2-\frac{1}{2\gamma_k}\|x-z_t\|^2\Big)+p_k\ell_{f}(\tilde x_{k},y_t)\nn \\
    &\quad\,-\frac{\alpha_k}{2\gamma_k}(1-L_{\nabla f}\alpha_k\gamma_k)\|z_{t-1}-z_t\|^2+\alpha_k\langle \delta_t, x-z_t\rangle
+ p_k\psi(\tilde x_{k})+\alpha_kC_5\epsilon. \label{l7-3}
\end{align}
In addition, observe that
\begin{align}
&p_k\ell_{f}(\tilde x_{k},y_t)-\frac{\alpha_k}{2\gamma_k}(1-L_{\nabla f}\alpha_k\gamma_k)\|z_{t-1}-z_t\|^2+\alpha_k\langle \delta_t, x-z_t\rangle \nn \\
 &=p_k\ell_{f}(\tilde x_{k},y_t)-\frac{\alpha_k}{2\gamma_k}(1-L_{\nabla f}\alpha_k\gamma_k)\|z_{t-1}-z_t\|^2+\alpha_k\langle \delta_t, z_{t-1}-z_t\rangle+\alpha_k\langle \delta_t, x-z_{t-1}\rangle \nn\\
    &\leq  p_k\ell_{f}(\tilde x_{k},y_t)+\frac{\alpha_k\gamma_k\|\delta_t\|^2}{2(1-L_{\nabla f}\alpha_k\gamma_k)}
    +\alpha_k\langle \delta_t, x-z_{t-1}\rangle, \label{l7-4}
\end{align}
where the last relation follows from \eqref{l7-1} and the Young's inequality.  Using \eqref{mean-var}, \eqref{l7-2}, the definition of $\ell_{f}(\cdot,\cdot)$, and the convexity of $f$, we obtain that
\begin{align}
&p_k\ell_{f}(\tilde x_{k},y_t)+\frac{\alpha_k\gamma_k\bbE[\|\delta_t\|^2|\Xi_{k,t-1}]}{2(1-L_{\nabla f}\alpha_k\gamma_k)}+\alpha_k\bbE[\langle \delta_t, x-z_{t-1}\rangle|\Xi_{k,t-1}]\nn\\
&\overset{\eqref{mean-var}}{\leq} p_k\ell_{f}(\tilde x_{k},y_t)+\frac{L_{\nabla f}\alpha_k\gamma_k}{1-L_{\nabla f}\alpha_k\gamma_k}\left(f(\tilde x_{k}) - \ell_{f}(\tilde x_{k},y_t)\right)\nn\\
&=\Big(p_k-\frac{L_{\nabla f}\alpha_k\gamma_k}{1-L_{\nabla f}\alpha_k\gamma_k}\Big)\ell_f(\tilde x_{k},y_t)+\frac{L_{\nabla f}\alpha_k\gamma_k}{1-L_{\nabla f}\alpha_k\gamma_k}f(\tilde x_{k})
\leq  p_kf(\tilde x_{k}),\label{l7-5}
\end{align}
where the last inequality follows from \eqref{l7-2} and the convexity of $f$. It then follows from \eqref{l7-4} and \eqref{l7-5} that
\begin{align*}
&\bbE\Big[p_k\ell_{f}(\tilde x_{k},y_t)-\frac{\alpha_k}{2\gamma_k}(1-L_{\nabla f}\alpha_k\gamma_k)\|z_{t-1}-z_t\|^2+\alpha_k\langle \delta_t, x-z_t\rangle\Big] \nn \\
    &\overset{\eqref{l7-4}}{\leq}  \bbE\Big[p_k\ell_{f}(\tilde x_{k},y_t)+\frac{\alpha_k\gamma_k\|\delta_t\|^2}{2(1-L_{\nabla f}\alpha_k\gamma_k)} +\alpha_k\langle \delta_t, x-z_{t-1}\rangle\Big] \overset{\eqref{l7-5}}{\leq} p_k\bbE[f
(\tilde x_{k})].
\end{align*}
Using this inequality and taking expectation on both sides of \eqref{l7-3}, we obtain that for all $1\leq t\leq T_k$, 
\[
 \bbE[F(x_t)] \leq(1-\alpha_k-p_k)\bbE[F(x_{t-1})]+\alpha_k\Big(F(x)+\frac{1}{2\gamma_k}\bbE[\|x-z_{t-1}\|^2]-\frac{1}{2\gamma_k}\bbE[\|x-z_t\|^2]\Big)+p_k\bbE[F(\tilde x_{k})]+\alpha_kC_5\epsilon.
\]
The relation \eqref{l7-exp} then follows from this inequality by rearranging the terms.
\end{proof}

We are now ready to prove Theorem \ref{thm:finite-sum}.

\begin{proof}[\textbf{Proof of Theorem \ref{thm:finite-sum}}]
We first show that the first two relations in \eqref{xk-prop} hold. Indeed, it follows from \eqref{omega-eps} and Lemma \ref{iter-property}(i) that $\dist(c(\tx_k),-\cK) \leq \epsilon$. Hence, there exists $u\in -\cK$ such that $\|c(\tx_k)-u\| \leq \epsilon$. Let $\lambda^*$ be the minimum-norm optimal Lagrange multiplier of problem \eqref{prob1}.  Then, one has that $\lambda^*\in\cK^*$ and 
\[
F^*= \min_x F(x)+\langle \lambda^*, c(x)\rangle.
\]
Using these facts, $u\in -\cK$, and $\|c(\tx_k)-u\| \leq \epsilon$, we obtain that
\[
\langle \lambda^*, c(\tilde x_k)\rangle \leq \langle \lambda^*, c(\tilde x_k)-u\rangle \leq  \|\lambda^*\| \epsilon, \qquad
F(\tx_k)+\langle \lambda^*, c(\tx_k)\rangle \geq F^*.
\]
In addition, by the definition of $D_\lambda$  in this theorem, one has $D_\lambda=\|\lambda^*\|$. Using these relations, we have
\[
F(\tx_k)  \geq F^*- \langle \lambda^*, c(\tx_k)\rangle  \geq F^* -  \|\lambda^*\| \epsilon = F^* -  D_\lambda \epsilon.
\]
Hence, the first two relations in \eqref{xk-prop} hold as desired.

We next show that the third relation in \eqref{xk-prop} holds. To this end, let $\tx^*\in\arg\min_{x\in\Omega_\epsilon} F(x)$. One can observe from  \eqref{prob1} and \eqref{omega-eps} that $F^*=\min_{x\in\Omega} F(x)$ and $\Omega \subseteq \Omega_\epsilon$. It then follows that 
$F(\tx^*)\leq F^*$.  In addition, one can verify that $\{q_i\}$, $\{\alpha_k\}$, $\{\gamma_k\}$, and $\{p_k\}$ that are given in \eqref{def-para} satisfy the assumptions stated in Lemma \ref{l-iter2}. Hence,  Lemma \ref{l-iter2} holds for such parameters.  It then follows from \eqref{l7-exp} with $x=\tx^*$ that for all $k \geq 1$, 
\begin{align}
    \frac{\gamma_k}{\alpha_k}\bbE[F(x_t)-F(\tx^*)]&\leq 
 \frac{\gamma_k}{\alpha_k}(1-\alpha_k-p_k)\bbE[F(x_{t-1})-F(\tx^*)]+\frac{\gamma_kp_k}{\alpha_k} \bbE[F(\tilde x_{k})-F(\tx^*)]\nn\\
 &\quad +\frac{1}{2}\big(\bbE[\|z_{t-1}-\tx^*\|^2]-\bbE[\|z_t-\tx^*\|^2]\big) +\gamma_kC_5\epsilon \qquad \forall 1 \leq t \leq T_k. \label{l9-1}
\end{align}
Summing up \eqref{l9-1} from $t=1$ to $T_k$, and using the fact that $x_0=\tilde x_{k}$, $z_0=\tilde z_{k}$, $z_{T_k}=\tilde z_{k+1}$, we have
\begin{align}
\sum_{t=1}^{T_k}\frac{\gamma_k}{\alpha_k}\bbE[F(x_t)-F(\tx^*)]
&\leq \frac{\gamma_k}{\alpha_k}(1-\alpha_k-p_k)\bbE[F(x_{0})-F(\tx^*)]+ \sum_{t=1}^{T_k-1}\frac{\gamma_k}{\alpha_k}(1-\alpha_k-p_k)\bbE[F(x_t)-F(\tx^*)]\nn\\
&\quad+\frac{T_k\gamma_kp_k}{\alpha_k}\bbE[F(\tilde x_{k})-F(\tx^*)]+\frac{1}{2}\big(\bbE[\|z_{0}-\tx^*\|^2]-\bbE[\|z_{T_k}-\tx^*\|^2]\big)+\gamma_kT_kC_5\epsilon\nn\\
&= \frac{\gamma_k}{\alpha_k}(1-\alpha_k-p_k)\bbE[F(\tilde x_{k})-F(\tx^*)]+ \sum_{t=1}^{T_k-1}\frac{\gamma_k}{\alpha_k}(1-\alpha_k-p_k)\bbE[F(x_{t})-F(\tx^*)]\nn\\
&\quad+\frac{T_k\gamma_kp_k}{\alpha_k}\bbE[F(\tilde x_{k})-F(\tx^*)]
+\frac{1}{2}\big(\bbE[\|\tilde z_{k}-\tx^*\|^2]-\bbE[\|\tilde z_{k+1}-\tx^*\|^2]\big)+\gamma_kT_kC_5\epsilon.\label{l9-sum}
\end{align}

Let $\cQ_k$ and $\mathcal{R}_k$ be defined as
\beq \label{Lk}
    \cQ_k :=\frac{\gamma_k}{\alpha_k}+(T_k-1)\frac{\gamma_k(\alpha_k+p_k)}{\alpha_k}, \qquad
    \mathcal{R}_k :=\frac{\gamma_k}{\alpha_k}(1-\alpha_k)+(T_k-1)\frac{\gamma_kp_k}{\alpha_k} \qquad \forall k \geq 1,
\eeq
where $\alpha_k$, $\gamma_k$, $p_k$, and $T_k$ are given in \eqref{def-para}. Notice from \eqref{Lk}, the definition of $\theta_t$ in \eqref{def-para}, and 
Algorithm \ref{alg2} that $\cQ_k=\sum_{t=1}^{T_k}\theta_t$ and $\tilde x_{k+1}=\sum_{t=1}^{T_k}(\theta_tx_t)/\sum_{t=1}^{T_k}\theta_t$. Using  these,  \eqref{def-para}, \eqref{l9-sum}, and the convexity of $F$, we have
    \begin{align}
        & \cQ_k\bbE[F(\tilde x_{k+1})-F(\tx^*)]=\Big(\sum_{t=1}^{T_k}\theta_t\Big)\bbE[F(\tilde x_{k+1})-F(\tx^*)]\leq\sum_{t=1}^{T_k}\theta_t\bbE[F(x_t)-F(\tx^*)]\nn\\
        &\overset{\eqref{def-para}}=\sum_{t=1}^{T_k}\frac{\gamma_k}{\alpha_k}\bbE[F(x_t)-F(\tx^*)]-\sum_{t=1}^{T_k-1}\frac{\gamma_k}{\alpha_k}(1-\alpha_k-p_k)\bbE[F(x_t)-F(\tx^*)]\nn \\
        &\overset{\eqref{l9-sum}}{\leq} \Big(\frac{\gamma_k}{\alpha_k}(1-\alpha_k)+(T_k-1)\frac{\gamma_kp_k}{\alpha_k}\Big)\bbE[F(\tilde x_{k})-F(\tx^*)]
    +\frac{1}{2}\big(\bbE[\|\tilde z_{k}-\tx^*\|^2]-\bbE[\|\tilde z_{k+1}-\tx^*\|^2]\big)+\gamma_kT_kC_5\epsilon \nn \\
& \overset{\eqref{Lk}}{=}\mathcal{R}_k\bbE[F(\tilde x_{k})-F(\tx^*)]
    +\frac{1}{2}\big(\bbE[\|\tilde z_{k}-\tx^*\|^2]-\bbE[\|\tilde z_{k+1}-\tx^*\|^2]\big)+\gamma_kT_kC_5\epsilon. \label{l9-2}
    \end{align}

We claim that $\cQ_k-\mathcal{R}_{k+1}\geq 0$ for all $k\geq1$. Indeed, notice from \eqref{def-para} that if $1\leq k< k_0$,  $\alpha_{k+1}=\alpha_k=1/2$, $p_k=1/2$, $\gamma_{k+1}=\gamma_k$, and $T_{k+1}=2T_k$. Using these relations and \eqref{Lk}, we have that for all $1\leq k< k_0$,
\begin{align*}
    &\cQ_k-\mathcal{R}_{k+1}
    \overset{\eqref{Lk}}{=}\frac{\gamma_k}{\alpha_k}+(T_k-1)\frac{\gamma_k(\alpha_k+p_k)}{\alpha_k}-\frac{\gamma_{k+1}}{\alpha_{k+1}}(1-\alpha_{k+1})-(T_{k+1}-1)\frac{\gamma_{k+1}p_{k+1}}{\alpha_{k+1}}\\
    &= \frac{\gamma_k}{\alpha_k}\Big(1+(T_k-1)(\alpha_k+p_k)-(1-\alpha_k)-(2T_k-1)p_k\Big)=0,
\end{align*}
and hence the above claim holds for all $1\leq k< k_0$. On the other hand, if $k\geq k_0$, notice from $k_0 \geq 1$ and  \eqref{def-para} that $T_k=T_{k+1}=T_{k_0}\geq1$, $\alpha_k=2/(k-k_0+4)$, $p_{k}=1/2$, and $\gamma_k=1/(3L_{\nabla f}\alpha_k)=(k-k_0+4)/(6L_{\nabla f})$. Using these relations and \eqref{Lk}, we obtain that for all $k\geq k_0$,
\begin{align*}
    &\cQ_k-\mathcal{R}_{k+1}\overset{\eqref{Lk}}{=}\frac{\gamma_k}{\alpha_k}+(T_{k}-1)\frac{\gamma_k(\alpha_k+p_k)}{\alpha_k}-\frac{\gamma_{k+1}}{\alpha_{k+1}}(1-\alpha_{k+1})-(T_{k+1}-1)\frac{\gamma_{k+1}p_{k+1}}{\alpha_{k+1}} \\
    &=\frac{\gamma_k}{\alpha_k}+(T_{k_0}-1)\frac{\gamma_k(\alpha_k+p_k)}{\alpha_k}-\frac{\gamma_{k+1}}{\alpha_{k+1}}(1-\alpha_{k+1})-(T_{k_0}-1)\frac{\gamma_{k+1}p_{k+1}}{\alpha_{k+1}} \\
    &=\frac{\gamma_k}{\alpha_k}-\frac{\gamma_{k+1}}{\alpha_{k+1}}(1-\alpha_{k+1})+(T_{k_0}-1)\Big(\frac{\gamma_k(\alpha_k+p_k)}{\alpha_k}-\frac{\gamma_{k+1}p_{k+1}}{\alpha_{k+1}}\Big) \\
    &=\frac{(k-k_0+4)^2}{12L_{\nabla f}}-\frac{(k-k_0+5)^2}{12L_{\nabla f}}\Big(1-\frac{2}{k-k_0+5}\Big) \\
    &\quad+(T_{k_0}-1)\Big[\frac{(k-k_0+4)^2}{12L_{\nabla f}}\Big(\frac{2}{k-k_0+4}+\frac{1}{2}\Big)-\frac{(k-k_0+5)^2}{24L_{\nabla f}}\Big] \\
    &=\frac{(k-k_0+4)^2-(k-k_0+5)^2+2(k-k_0+5)}{12L_{\nabla f}} \\
    &\quad\,+\frac{(T_{k_0}-1)\big[4(k-k_0+4)+(k-k_0+4)^2-(k-k_0+5)^2\big]}{24L_{\nabla f}} \\
    &=\frac{1}{12L_{\nabla f}}+\frac{(T_{k_0}-1)(2(k-k_0)+7)}{24L_{\nabla f}}\geq 0,
\end{align*}
where the second equality is due to $T_k=T_{k+1}=T_{k_0}$, and the fourth equality follows from $\alpha_k=2/(k-k_0+4)$, $p_{k}=1/2$, and $\gamma_k=(k-k_0+4)/(6L_{\nabla f})$. Hence, the above claim also holds for all $k \geq k_0$.

Recall from Lemma \ref{iter-property}(i) that $\tilde x_k\in\Omega_\epsilon$ for all $k\geq1$. It follows from this fact and the definition $\tx^*$ that $F(\tilde x_k)-F(\tx^*)\geq 0$. Also, note from \eqref{def-para} and \eqref{Lk} that $\mathcal{R}_1=2/(3L_{\nabla f})$.  Summing up \eqref{l9-2} with $k$ replaced by $i$ from $i=1$ to $k-1$, and using the fact that $F(\tilde x_i)-F(\tx^*)\geq 0$, $\cQ_i-\mathcal{R}_{i+1}\geq 0$ for all $i\geq 1$, $\tilde x_1=\tilde z_1=x^0$, and $\mathcal{R}_1=2/(3L_{\nabla f})$, we have
\begin{align}
\cQ_{k-1}\bbE[F(\tilde x_k)-F(\tx^*)]&
\leq \sum_{i=1}^{k-2}(\mathcal{R}_{i+1}-\cQ_i)\bbE[F(\tilde x_{i+1})-F(\tx^*)]+\mathcal{R}_1(F(\tilde x_1)-F(\tx^*))\nn\\
&\quad\ +\frac{1}{2}\|\tilde z_1-\tx^*\|^2+C_5\epsilon\sum_{i=1}^{k-1}\gamma_iT_i \nn\\
&\leq\frac{2}{3L_{\nabla f}}(F(x^0)-F(\tx^*))+\frac{1}{2}\|x^0-\tx^*\|^2+C_5\epsilon\sum_{i=1}^{k-1}\gamma_iT_i.    \label{l9-4}
\end{align}
In addition, recall that $k_0=\lfloor \log_2 s\rfloor+1$ and $T_{k_0}=2^{k_0-1}$, which imply that  $s/2\leq T_{k_0} \leq s$. By this relation,  \eqref{def-para}, and \eqref{Lk}, one has that for all $k>k_0$,
\begin{align}
    \cQ_{k-1}&\overset{\eqref{Lk}}{=}\frac{\gamma_{k-1}}{\alpha_{k-1}}\big(1+(T_{k-1}-1)(\alpha_{k-1}+p_{k-1})\big) \nn \\
    &\overset{\eqref{def-para}}{=}\frac{1}{3L_{\nabla f}\alpha_{k-1}^2}\left(1+(T_{k-1}-1)\Big(\alpha_{k-1}+\frac{1}{2}\Big)\right)
    =\frac{T_{k-1}-1}{3L_{\nabla f}\alpha_{k-1}}+\frac{T_{k-1}+1}{6L_{\nabla f}\alpha_{k-1}^2} \nn \\
    &\overset{\eqref{def-para}}{=}\frac{(T_{k_0}-1)(k-k_0+3)}{6L_{\nabla f}}+\frac{(T_{k_0}+1)(k-k_0+3)^2}{24L_{\nabla f}}\nn \\
    &\geq\frac{(s-2)(k-k_0+3)}{12L_{\nabla f}}+\frac{(s+2)(k-k_0+3)^2}{48L_{\nabla f}} \nn \\
    & =\frac{(s-2)(k-k_0+3)}{12L_{\nabla f}}+\frac{2(k-k_0+3)^2}{48L_{\nabla f}} + \frac{s(k-k_0+3)^2}{48L_{\nabla f}} \nn \\
    &\geq \frac{(s-2)(k-k_0+3)}{12L_{\nabla f}}+\frac{8(k-k_0+3)}{48L_{\nabla f}} + \frac{s(k-k_0+3)^2}{48L_{\nabla f}} \nn \\
    & = \frac{s(k-k_0+3)}{12L_{\nabla f}}+ \frac{s(k-k_0+3)^2}{48L_{\nabla f}} 
     \geq \frac{s(k-k_0+3)^2}{48L_{\nabla f}}, \label{Q-bound1} \\
    \sum_{i=1}^{k-1}\gamma_iT_i&=\sum_{i=1}^{k_0-1}\gamma_iT_i+\sum_{i=k_0}^{k-1}\gamma_iT_i 
    \overset{\eqref{def-para}}{=} \frac{2}{3L_{\nabla f}}\sum_{i=1}^{k_0-1}2^{i-1}+\frac{T_{k_0}}{3L_{\nabla f}}\sum_{i=k_0}^{k-1}\frac{i-k_0+4}{2} \nn \\
    &\leq \frac{2^{k_0}}{3L_{\nabla f}}+\frac{T_{k_0}}{3L_{\nabla f}}\frac{(k-k_0+7)(k-k_0)}{4}
    = \frac{T_{k_0}\big((k-k_0)^2+7(k-k_0)+8\big)}{12L_{\nabla f}} \nn \\
    &\leq \frac{T_{k_0}\big(2(k-k_0)^2+12(k-k_0)+18\big)}{12L_{\nabla f}}\leq\frac{s(k-k_0+3)^2}{6L_{\nabla f}}. \label{Q-bound2}
\end{align}
Notice that $x^0, \tx^*\in X$. It then follows from Assumption \ref{a1} that $F(x^0)-F(\tx^*) \leq D_F$ and  
$\|x^0-\tx^*\| \leq D_X$. Using these relations, \eqref{l9-4}, \eqref{Q-bound1}, and  \eqref{Q-bound2}, we obtain that
\begin{align*}
\bbE[F(\tilde x_k)-F(\tx^*)]&\leq \frac{16}{s(k-k_0+3)^2}\Big(2(F(x^0)-F(\tx^*))+\frac{3L_{\nabla f}}{2}\|x^0-\tx^*\|^2\Big)+8C_5\epsilon\nn\\
&\leq \frac{16}{s(k-k_0+3)^2}\Big(2D_F + \frac{3L_{\nabla f}  D_X^2}{2}\Big)+8C_5\epsilon,
\end{align*}
which together with \eqref{c45} and $F(\tx^*)\leq F^*$  implies that the third relation in \eqref{xk-prop} holds.

In the remainder of this proof, we show that \eqref{N-bound} holds. To this end, let $N_{k,f}$ denote the number of evaluations of $\nabla f_i$ performed in the $k$th outer iteration of Algorithm \ref{alg2} for any $k \geq 1$, and let $N_{k,t}$ denote the number of evaluations of $\nabla c$ performed in the $t$-th inner iteration at the $k$th outer iteration of Algorithm \ref{alg2} for all $1 \leq t \leq T_k$ and $k \geq 1$.

Observe from Algorithm \ref{alg2} that $N_{k,f}=s+T_k$ for all $k\geq 1$.  We next bound $N_{k,t}$ for all $1 \leq t \leq T_k$ and $k\geq 1$.  To this end,   let $(x^*_{k,t},\lambda^*_{k,t})$ be a pair of optimal solutions of problem \eqref{subproblem} and its dual, and define
\begin{align}
&  D_{k,t}= \max\big\{\big\|(z_{t-1},\lambda_{t-1} \min\{R/\|\lambda_{t-1}\|, 1\})-(x^*_{k,t},\lambda^*_{k,t})\big\|, 1\big\}, \label{Dkt1} \\ 
& B_{k,t} = 1/\gamma_k + L_{\nabla c}\|\lambda^*_{k,t}\| + L_{\nabla c}D_{k,t},  \label{Bkt} \\
&\bar{C}_{k,t} = C \rho_0+B_{k,t}+L_{\nabla c}\rho_0\eta_0/(1-\zeta\sigma)+1. \label{Ckt1} 
\end{align}
By \eqref{DX}, \eqref{lambda*}, \eqref{Dkt1}, and $z_{t-1}, x^*_{k,t}\in X$, one has 
\beq \label{Dkt2}
 D_{k,t}  \leq \max\big\{\|z_{t-1}-x^*_{k,t}\| +\|\lambda_{t-1} \min\{R/\|\lambda_{t-1}\|, 1\}\|+\|\lambda^*_{k,t}\|, 1\big\} \leq D_X + R + \Lambda.
\eeq
Observe from \eqref{def-para} that $1/\gamma_k=3L_{\nabla f}\alpha_k\leq 3L_{\nabla f}/2$. In view of this,  \eqref{lambda*},  \eqref{Bkt}, and  \eqref{Dkt2}, we obtain that
\[
B_{k,t} \leq  3L_{\nabla f}/2  + L_{\nabla c}\Lambda+ L_{\nabla c}(D_X + R + \Lambda)= 3L_{\nabla f}/2  + L_{\nabla c}(D_X + R +2 \Lambda).
\]
Using this relation, \eqref{C2}, \eqref{Ckt1}, and the definition of $C$ in \eqref{DX}, we have
\begin{align}
 \bar{C}_{k,t} & \overset{\eqref{DX}\eqref{Ckt1}}{=} (L_{\nabla c}C_c+L_c^2)\rho_0+B_{k,t}+L_{\nabla c}\rho_0\eta_0/(1-\zeta\sigma)+1 \overset{\eqref{C2}}{\leq} C_2. \label{Ckt2}
\end{align}
Note from \eqref{subproblem} that the smooth component of $F_{k,t}(x)$ is $\|x-z_{t-1}+\gamma_kg_t\|^2/(2\gamma_k)$, which is $1/\gamma_k$-smooth and $1/\gamma_k$-strongly convex in $\rr^n$.  Moreover, notice from Algorithm \ref{alg2} that Algorithm \ref{finitesum:alg1} with $\bepsilon = \epsilon$, $\mu = 1/\gamma_k$, $(x_0,\lambda_0) =(z_{t-1},\lambda_{t-1} \min\{R/\|\lambda_{t-1}\|, 1\})$ is called to solve problem \eqref{subproblem}. By these observations, \eqref{C2}, \eqref{C3}, \eqref{Dkt2}, \eqref{Ckt2}, and Theorem \ref{thm:sub-complexity}, one has
\begin{align}
N_{k,t} &\overset{\text{Theorem } \ref{thm:sub-complexity}}{\leq} \frac{16\zeta\sqrt{\bar{C}_{k,t}}}{\sqrt{1/\gamma_k}(\sqrt{\zeta}-1)}\log\frac{2\zeta  \bar{C}_{k,t} D_X}{\eta_0\sigma}\left\lceil\max\left\{1,  \log_\zeta\frac{2(D_{k,t}+\rho_0\eta_0)}{(1-\zeta\sigma)\epsilon}\right\}\right\rceil\max\left\{1,  \sqrt{\frac{2(D_{k,t}+\rho_0\eta_0)}{(1-\zeta\sigma)\epsilon}}\right\} \nn \\
&\leq  L_{\nabla f}^{\frac{1}{2}}C_3 \sqrt{\gamma_k} \left\lceil \max\left\{1, \log_\zeta\frac{C_4}{\epsilon}\right\}\right\rceil \max\left\{1,  \sqrt{\frac{C_4}{\epsilon}} \right\}. \label{Nkt}
\end{align}

In addition, since $k_0=\lfloor \log_2 s \rfloor+1$, it follows that $k_0 \geq 1$, and hence $i-k_0+4 \leq 4i$ for all $i\geq 1$. Using this and the definitions of $\alpha_i$ and $\gamma_i$ in \eqref{def-para}, we have
\begin{align}
 \sum_{i=1}^{k}\sqrt{\gamma_i} &=(3 L_{\nabla f})^{-\frac{1}{2}}\sum_{i=1}^{k}\alpha_i^{-\frac{1}{2}}=(3 L_{\nabla f})^{-\frac{1}{2}}\Big(\sum_{i=1}^{k_0} \sqrt{2}+\sum_{i=k_0+1}^{k} \sqrt{(i-k_0+4)/2}\,\Big) 
\leq (3 L_{\nabla f})^{-\frac{1}{2}}\sum_{i=1}^{k}\sqrt{2i} \nn \\ 
&\leq \frac{2}{3}(3 L_{\nabla f}/2)^{-\frac{1}{2}}(k+1)^{\frac{3}{2}} \leq  \frac{2}{3}(3 L_{\nabla f}/2)^{-\frac{1}{2}}(2k)^{\frac{3}{2}} \leq 2 L_{\nabla f}^{-\frac{1}{2}} 
k^{\frac{3}{2}}, \label{sum-gamma}
\end{align}
where the last inequality is due to the fact that $\sum_{i=1}^{k}\sqrt{i}\leq\int_{\tau=0}^{k}\sqrt{\tau+1}d\tau <2 (k+1)^{3/2}/3$. 

Let $N_0$ denote the number of evaluations of $\nabla c$ performed by Algorithm \ref{alg3} applied to the problem $\min_{x \in X} \dist^2(c(x),-\cK)$ for finding $x^0\in X$ satisfying $\dist(c(x^0),-\cK)\leq\epsilon$. 
One can verify that the function $\dist^2(c(x),-\cK)$ is $2(L_c^2+L_{\nabla c}C_c)D_X$-smooth on $X$. 
In addition, notice that one evaluation of $\nabla c$ is performed by Algorithm \ref{alg3} per iteration. It then follows from \eqref{C2} and Theorem \ref{thm:cplx-fom-coro} with $L_{\nabla \tf}=2(L_c^2+L_{\nabla c}C_c)D_X$ and $\varepsilon=\epsilon^2$ that 
\beq \label{N0}
N_0 \leq 2 D_X^{3/2} \sqrt{L_c^2+L_{\nabla c}C_c}\epsilon^{-1} \overset{\eqref{C2}}{=} C_1 \epsilon^{-1}.
\eeq

Recall that $N_f$ and $N_c$ denote the number of evaluations of $\nabla f_i$ and $\nabla c$  during the $k$ outer iterations of Algorithm \ref{alg2} for all $k \geq 1$, respectively. Notice from \eqref{def-para} that $T_i\leq T_{k_0}\leq s$ for all $i \geq 1$. Using this, \eqref{Nkt}, \eqref{sum-gamma}, \eqref{N0}, and Algorithm \ref{alg2}, we can observe that for all $k \geq 1$,
\begin{align*}
N_f &= \sum_{i=1}^{k}(s+T_i)\leq2sk, \\
 N_c&=N_0+\sum_{i=1}^{k}\sum^{T_i}_{t=1}N_{i,t} \overset{\eqref{Nkt}}{\leq} N_0+L_{\nabla f}^{\frac{1}{2}}C_3  \left\lceil \max\left\{1, \log_\zeta\frac{C_4}{\epsilon}\right\}\right\rceil \max\left\{1,  \sqrt{\frac{C_4}{\epsilon}} \right\} \sum_{i=1}^{k} T_i   \sqrt{\gamma_i}  \\
&\leq N_0+s L_{\nabla f}^{\frac{1}{2}}C_3  \left\lceil \max\left\{1, \log_\zeta\frac{C_4}{\epsilon}\right\}\right\rceil \max\left\{1,  \sqrt{\frac{C_4}{\epsilon}} \right\} \sum_{i=1}^{k}  \sqrt{\gamma_i},
\end{align*}
which together with \eqref{sum-gamma} and \eqref{N0} implies that \eqref{N-bound} holds.
\end{proof}

\section{Concluding remarks}\label{sec:conclude}

In this paper, we propose an inexact stochastic variance-reduced gradient method for solving a class of constrained convex finite-sum optimization problems. The method computes an $\epsilon$-\textit{surely feasible stochastic optimal} ($\epsilon$-SFSO) solution using $\cO\big(s \log s + s^{1/2} \epsilon^{-1/2}\big)$ and $\widetilde{\cO}\big(s (\log s)^{3/2}\epsilon^{-1/2} + s^{1/4} \epsilon^{-5/4}\big)$ evaluations of $\nabla f_i$ and $\nabla c$, respectively. These bounds significantly improve upon the previously best-known first-order oracle complexity for obtaining an $\epsilon$-SFSO solution. 

It is worth noting that the proposed method can be extended to the case where the objective function is additionally strongly convex. In this setting, the extended method can compute an $\epsilon$-SFSO solution using
\[
\cO\big(s\log s + s^{1/2}\log \epsilon^{-1}\big)
\quad \text{and} \quad
\cO\Big(\epsilon^{-1} + s^{1/2}\epsilon^{-1/2} (\log \epsilon^{-1})^2 + s \log s\, \epsilon^{-1/2} \log \epsilon^{-1} \Big)
\]
evaluations of $\nabla f_i$ and $\nabla c$, respectively. These bounds further improve upon those obtained under the merely convex setting.

\appendix

\section{Accelerated proximal gradient methods for unconstrained convex optimization} \label{appendix}

In this part, we review accelerated proximal gradient (APG) methods for solving a class of unconstrained convex optimization problems of the form
\begin{equation}\label{eq:uncons-cvx-opt}
\tF^* = \min_{x\in\rr^n} \left\{\tF(x) := \tf(x) + \psi(x) \right\}, 
\end{equation}
where $\psi$ is a closed convex function with compact domain denoted by $X$, and $\tf$ is an $L_{\nabla \tf}$-smooth convex function on $X$. Throughout this section, we define $D_X=\max\{1,\max_{x,y\in X}\|x-y\|\}$.

We first present an APG method for solving \eqref{eq:uncons-cvx-opt} in which $\tf$ is convex but not necessarily strongly convex. This method is a variant of Nesterov’s optimal first-order  methods \cite{Nesterov13,Nesterov05}, and has been studied in, for example, \cite[Section~3]{Tseng08}.

\begin{algorithm}[H]
\caption{An APG method  for  problem \eqref{eq:uncons-cvx-opt} with general convex $\tf$}\label{alg3}
\begin{algorithmic}[1]
\REQUIRE  $x_{0} \in X$. \\
\STATE Set $z_0=x_0$.
\FOR{$k = 0,1,\dots$} 
    \STATE Compute
    \[
    \begin{aligned}
   y_k &= (k x_k + 2z_k)/(k+2), \\
   z_{k+1} &= \arg\min_{x} \Big\{ \langle \nabla \tf(y_k), x \rangle + \psi(x) + \frac{L_{\nabla \tf}}{k+2}\|x - z_k\|^2\Big\}, \\
   x_{k+1} &= (k x_k + 2z_{k+1})/(k+2).  
    \end{aligned}
    \]
\ENDFOR
\end{algorithmic}
\end{algorithm}

The following result provides the iteration complexity of Algorithm~\ref{alg3} for finding an $\varepsilon$-optimal solution to problem~\eqref{eq:uncons-cvx-opt}. It follows directly from \cite[Corollary~1]{Tseng08}, and thus its proof is omitted.

\begin{thm}
	\label{thm:cplx-fom-coro}
	For any given $\varepsilon>0$, Algorithm \ref{alg3} finds an approximate solution $x_k$ of problem \eqref{eq:uncons-cvx-opt} such that $\tF(x_k) - \tF^*\leq  \varepsilon$ in no more than $\ceil{D_X\sqrt{2L_{\nabla \tf}/\varepsilon}\,}$ iterations.
\end{thm}

We next present an optimal APG method \cite[Algorithm 4]{lu2023iteration}, which is a slight variant of Nesterov's optimal first-order methods \cite{LLX15,Nesterov13} for solving strongly convex problem \eqref{eq:uncons-cvx-opt}, where $\tf$ is further assumed to be $\mu_\tf$-strongly convex on $X$ with $0<\mu_\tf \leq L_{\nabla \tf}$.

\begin{algorithm}[H]
\caption{An APG method for problem \eqref{eq:uncons-cvx-opt} with strongly convex $\tf$}\label{alg4}
\begin{algorithmic}[1]
\REQUIRE $\varepsilon > 0$ and $\tx_{0} \in X$. \\
\STATE Set $x_0=\tx_0$ and compute 
\begin{equation*}
   x_{1} = \argmin_{x} \Big\{\langle\nabla \tf(x_{0}),x\rangle + \psi(x)+\frac{L_{\nabla \tf}}{2}\|x-x_{0}\|^2\Big\}.
\end{equation*}
\STATE Set $z_{1}=x_{1}$ and $\alpha=\sqrt{\mu_\tf/L_{\nabla \tf}}$.
\FOR{$k = 1,2,\dots$} 
    \STATE Compute
    \[
    \begin{aligned}
   y_k &= (x_k +\alpha z_k ) / (1 +\alpha ), \\
   z_{k+1} &= \arg\min_{x} \Big\{ \langle \nabla \tf(y_k), x \rangle + \psi(x) + \frac{\alpha L_{\nabla \tf}  }{2} \|x -\alpha y_k - (1 -\alpha)z_k \|^2 \Big\}, \\
   x_{k+1} &= (1 -\alpha)x_k +\alpha  z_{k+1}.
    \end{aligned}
    \]
     \STATE Compute 
        \begin{equation*}
            \hat{x}_{k+1}=\arg\min_{x} \Big\{ \langle \nabla \tf(x_{k+1}), x \rangle +  \psi(x) + \frac{L_{\nabla \tf}}{2} \|x -x_{k+1}\|^2 \Big\}.
        \end{equation*}
        \IF{$\|x_{k+1} - \hat{x}_{k+1}\| \leq \varepsilon/(2L_{\nabla \tf})$}
            \STATE Terminate the algorithm and output $\hat{x}_{k+1}$.
    \ENDIF
\ENDFOR
\end{algorithmic}
\end{algorithm}

The following result provides the iteration complexity of Algorithm \ref{alg4} for computing an approximately optimal solution to problem~\eqref{eq:uncons-cvx-opt}. Its proof follows from \cite[Theorem~1]{LLX15}, combined with standard arguments, and is therefore omitted.

\begin{thm} \label{thm:cplx-fom-str-cvx}
For any given $\varepsilon>0$, Algorithm~\ref{alg4} terminates in at most $\widetilde K$ iterations, and outputs an approximate solution $\hat{x}_{k+1}$ of problem \eqref{eq:uncons-cvx-opt} satisfying $\dist(0, \partial \tF(\hat{x}_{k+1})) \leq \varepsilon$, where
\[
\widetilde K = \left\lceil \sqrt{\frac{L_{\nabla \tf}}{\mu_\tf}}\ \right\rceil \max\left\{ 1, \left\lceil 2\log\frac{2L_{\nabla \tf}D_X}{\varepsilon}\right\rceil\right\}.
\]
\end{thm}

\end{document}